\documentclass[onefignum,onetabnum,dvipsnames]{siamonline190516}
\usepackage{csquotes}
\usepackage{bibentry}
\usepackage{cite}
\usepackage{epigraph}

\usepackage{bm}
\usepackage[mathscr]{euscript}
\usepackage{stmaryrd}
\SetSymbolFont{stmry}{bold}{U}{stmry}{m}{n} 
\usepackage[normalem]{ulem}
\usepackage{verbatim}
\usepackage{hyphenat} 
\usepackage{ulem}
\usepackage{cancel}

\usepackage{listings}
\usepackage{hyperref}
\usepackage{cleveref}

\crefname{lstlisting}{tensor}{listings}
\Crefname{lstlisting}{Tensor}{Listings}

\usepackage[space=true]{accsupp}
\newcommand{\copyablespace}{\BeginAccSupp{method=hex,unicode,ActualText=00A0}\
\EndAccSupp{}}

\usepackage{caption}
\usepackage{float}
\usepackage{graphicx}
\graphicspath{{./images/}}
\usepackage{subcaption}
\usepackage{tikz}
\usepackage{wrapfig}
\usepackage{pgfplots}
\pgfplotsset{compat=1.16}
\pgfplotsset{select coords between index/.style 2 args={ x filter/.code={
			\ifnum\coordindex<#1\fi
			\ifnum\coordindex>#2\fi } }}
\usepackage{booktabs}
\usepackage{ltablex}
\usepackage{multirow}
\usepackage{tabularx}
\usepackage{threeparttable}

\usepackage{adjustbox}

\usepackage{amsfonts,amsmath,amssymb,amsbsy}
\usepackage{mathtools}
\usepackage{etoolbox}

\usepackage{algorithm}
\usepackage[noend]{algpseudocode}
\Crefname{ALC@unique}{Line}{Lines} 
\algnewcommand{\LeftComment}[1]{\{#1\}\hfill}

\usepackage{enumitem}
\usepackage{siunitx}
\usepackage{url}

\usepackage{xparse}

\newcommand{\TODO}[2]{}
\newcommand{\EnableTODO}{%
	\renewcommand{\TODO}[2]{\textcolor{blue}{\textbf{\textit{[##1---##2]}}}}}
 \EnableTODO

\newcommand{\CHANGED}[1]{}
\newcommand{\EnableCHANGED}{%
	\renewcommand{\CHANGED}[1]{%
		\textcolor{Black}{##1}%
}}
 \EnableCHANGED

\NewDocumentCommand{\N}{o}{\mathbb{N}} \NewDocumentCommand{\R}{o}{\mathbb{R}}
\NewDocumentCommand{\F}{o}{\mathbb{F}} \NewDocumentCommand{\Z}{o}{\mathbb{Z}}
\NewDocumentCommand{\BigO}{o}{\mathcal{O}}
\NewDocumentCommand{\I}{o}{\mathcal{I}}

\makeatletter
\newcommand{\@giventhatstar}[2]{\left(#1\;\middle|\;#2\right)}
\newcommand{\@giventhatnostar}[3][]{#1(#2\;#1|\;#3#1)}
\newcommand{\giventhat}{\@ifstar\@giventhatstar\@giventhatnostar}
\makeatother

\NewDocumentCommand{\tns}{m}{\boldsymbol{\mathscr{#1}}}
\NewDocumentCommand{\mind}{o}{\boldsymbol{\mathbf{i}}}
\NewDocumentCommand{\MU}{o}{CP-APR-MU}
\NewDocumentCommand{\PDNR}{o}{CP-APR-PDNR}
\NewDocumentCommand{\PQNR}{o}{CP-APR-PQNR}
\NewDocumentCommand{\fsampnzs}{o}{p_{\tilde{\mat{F}}}}
\NewDocumentCommand{\fsampzs}{o}{q_{\tilde{\mat{F}}}}
\NewDocumentCommand{\gsampnzs}{o}{p_{\tilde{\mat{G}}}}
\NewDocumentCommand{\gsampzs}{o}{q_{\tilde{\mat{G}}}}

\NewDocumentCommand{\matbold}{m}{%
  \boldsymbol{\mathbf{#1}}%
} \NewDocumentCommand{\mat}{moo}{%
  \IfValueTF{#2}%
  {%
    \IfValueTF{#3}%
    {%
      \matbold{#1}(#2,#3)%
    }%
    {%
      \matbold{#1}_{#2}%
    }%
  }%
  {%
    \matbold{#1}%
  }%
}

\NewDocumentCommand{\vecbold}{m}{%
  \boldsymbol{\mathbf{#1}} } \RenewDocumentCommand{\vec}{moo}{%
  \IfValueTF{#2}%
  {
    \IfValueTF{#3}%
    {
      \vecbold{#1}_{#2}(\colon,\,#3)%
    }%
    {
      \vecbold{#1}_{#2}%
    }%
  }%
  {
    \vecbold{#1}%
  }%
}

\NewDocumentCommand{\ktns}{m}{\llbracket #1\rrbracket}

\NewDocumentCommand\xDeclarePairedDelimiter{mmm} {%
  \NewDocumentCommand#1{som}{%
    \IfNoValueTF{##2} {\IfBooleanTF{##1}{#2##3#3}{\mleft#2##3\mright#3}}
    {\mathopen{##2#2}##3\mathclose{##2#3}}%
  }%
} \DeclareMathOperator{\diag}{diag}       
\DeclareMathOperator{\FMS}{FMS}         
\DeclareMathOperator{\nnz}{nnz}         
\DeclareMathOperator{\nz}{nz}           
\DeclareMathOperator{\rank}{rank}       

\DeclarePairedDelimiter\abs{\lvert}{\rvert}%
\DeclarePairedDelimiter\norm{\lVert}{\rVert}%

\makeatletter
\let\oldabs\abs \def\abs{\@ifstar{\oldabs}{\oldabs*}}
\let\oldnorm\norm \def\norm{\@ifstar{\oldnorm}{\oldnorm*}}
\makeatother

\NewDocumentCommand{\eqwhere}{o}{\quad \text{where} \quad}
\NewDocumentCommand{\eqst}{o}{\quad \text{s.t.} \quad}
\NewDocumentCommand{\eqand}{o}{\quad \text{and} \quad}
\NewDocumentCommand{\eqfor}{o}{\quad \text{for} \quad}
\NewDocumentCommand{\eqforall}{o}{\quad \text{for all} \quad}
\NewDocumentCommand{\eqwith}{o}{\quad \text{with} \quad}
\NewDocumentCommand{\eqif}{o}{\quad \text{if} \quad}
\NewDocumentCommand{\eqow}{o}{\quad \text{otherwise} \quad}

\algblockdefx[IfOnce]{IfOnce}{EndIfOnce}
[1]{\textbf{if}~#1~\textbf{do only once}}[0]{\textbf{end if}}
\algnewcommand\CommentLine[1]{// #1}

\makeatletter
\let\save@mathaccent\mathaccent \newcommand*\if@single[3]{%
\setbox0\hbox{${\mathaccent"0362{#1}}^H$}%
\setbox2\hbox{${\mathaccent"0362{\kern0pt#1}}^H$}%
\ifdim\ht0=\ht2 #3\else #2\fi }
\newcommand*\rel@kern[1]{\kern#1\dimexpr\macc@kerna}
\newcommand*\widebar[1]{\@ifnextchar^{{\wide@bar{#1}{0}}}{\wide@bar{#1}{1}}}
\newcommand*\wide@bar[2]{\if@single{#1}{\wide@bar@{#1}{#2}{1}}{\wide@bar@{#1}{#2}{2}}}
\newcommand*\wide@bar@[3]{%
  \begingroup
  \def\mathaccent##1##2{%
    \let\mathaccent\save@mathaccent
    \if#32 \let\macc@nucleus\first@char \fi
    \setbox\z@\hbox{$\macc@style{\macc@nucleus}_{}$}%
    \setbox\tw@\hbox{$\macc@style{\macc@nucleus}{}_{}$}%
    \dimen@\wd\tw@ \advance\dimen@-\wd\z@
    \divide\dimen@ 3 \@tempdima\wd\tw@ \advance\@tempdima-\scriptspace
    \divide\@tempdima 10 \advance\dimen@-\@tempdima
    \ifdim\dimen@>\z@ \dimen@0pt\fi
    \rel@kern{0.6}\kern-\dimen@ \if#31
    \overline{\rel@kern{-0.6}\kern\dimen@\macc@nucleus\rel@kern{0.4}\kern\dimen@}%
      \advance\dimen@0.4\dimexpr\macc@kerna
      \let\final@kern#2%
      \ifdim\dimen@<\z@ \let\final@kern1\fi \if\final@kern1 \kern-\dimen@\fi
    \else
      \overline{\rel@kern{-0.6}\kern\dimen@#1}%
    \fi
  }%
  \macc@depth\@ne \let\math@bgroup\@empty \let\math@egroup\macc@set@skewchar
  \mathsurround\z@ \frozen@everymath{\mathgroup\macc@group\relax}%
  \macc@set@skewchar\relax \let\mathaccentV\macc@nested@a
  \if#31 \macc@nested@a\relax111{#1}%
  \else
    \def\gobble@till@marker##1\endmarker{}%
    \futurelet\first@char\gobble@till@marker#1\endmarker
    \ifcat\noexpand\first@char A\else \def\first@char{}%
    \fi
    \macc@nested@a\relax111{\first@char}%
  \fi
  \endgroup
}
\makeatother

\headers{Tensor Decompositions for Count Data}{J. M. Myers and D. M. Dunlavy}

\title{Tensor Decompositions for Count Data that Leverage Stochastic and
Deterministic Optimization}

\author{Jeremy M. Myers\thanks{Sandia National Laboratories
		(\email{jermyer@sandia.gov}, \email{dmdunla@sandia.gov})} \and Daniel M.
		Dunlavy\footnotemark[1] }

\ifpdf
\hypersetup{ pdftitle={Tensor Decompositions for Count Data that Leverage
	Stochastic and Deterministic Optimization}, pdfauthor={J. M. Myers and D. M.
	Dunlavy} }
\fi

\begin{document}

\maketitle

\begin{abstract}

There is growing interest to extend low-rank matrix decompositions to multi-way
arrays, or \emph{tensors}. One fundamental low-rank tensor decomposition is the
\emph{canonical polyadic decomposition (CPD)}. The challenge of fitting a
low-rank, nonnegative CPD model to Poisson-distributed count data is of
particular interest.
\CHANGED{%
Several popular algorithms use local search methods to approximate the
maximum likelihood estimator (MLE) of the Poisson CPD model.
This work presents two new algorithms that extend state-of-the-art
local methods for Poisson CPD. Hybrid GCP-CPAPR
combines
Generalized
Canonical Decomposition (GCP) with stochastic optimization
and CP Alternating Poisson Regression (CPAPR), a deterministic algorithm, to
increase the probability of converging to the MLE over either method used alone.
Restarted CPAPR
with \textsc{SVDrop} uses a heuristic based on the singular values of
the CPD 
model unfoldings to identify
convergence toward
optimizers that are not the MLE and restarts within
the feasible domain of the optimization problem, thus reducing overall computational
cost when using a multi-start strategy.
We provide empirical evidence
that indicates our approaches outperform existing methods with respect
to converging to the Poisson CPD MLE.%
}

\end{abstract}

\begin{keywords}
	tensor, canonical polyadic decomposition, GCP, CPAPR, count data, Poisson
\end{keywords}

\begin{AMS}
	15A69, 65F55
\end{AMS}

\capstartfalse
\begin{figure}[b!]
\centering
\includegraphics[width=\textwidth]{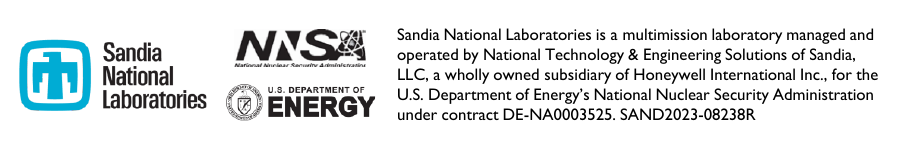}
\end{figure}
\capstarttrue

\section{Introduction}
\label{sec:intro}

Low-rank tensor decompositions in general, and the canonical polyadic
decomposition (CPD) specifically, \CHANGED{are important}
for multi-way data analysis\cite{Kolda09TensorDecompositionsApplications}. 
Fitting the parameters of a low-rank, nonnegative CPD
model to count data is often formulated as a nonlinear, nonconvex global
optimization problem. When the data are assumed to be Poisson-distributed, one
approach is to determine the optimal Poisson parameters that maximize the
likelihood of the data via tensor maximum likelihood
estimation\cite{Chi12TensorsSparsityNonnegative}. The global
optimizer to the optimization problem is the maximum likelihood estimator (MLE).
Since global optimization algorithms are often prohibitively expensive for
tensor data, great emphasis has been placed on developing efficient local
methods for finding the Poisson CPD parameters. In practice, local methods for
solving global optimization problems are often orchestrated in a multi-start
strategy---i.e., computing a set of approximations from many random starting
points---to increase the probability that the model best approximating the MLE
has been found. However, this approach demands significant computational
resources when high-confidence solutions are required and may lead to excessive
computations even for small problems. To mitigate this issue, we examine the
role of randomization and determinism in Poisson CPD solvers.

Our contributions are:
\begin{itemize}
	\item two Poisson CPD methods that compute the MLE with higher probability
	than current effective local search methods, and
	\item validation of our methods with open-source software on synthetic data.
\end{itemize}

\subsection{Hybrid GCP-CPAPR \textsc{(HybridGC)}}
Our first method is \textsc{HybridGC}, which uses a two-stage hybrid strategy
built from effective local methods.
\CHANGED{%
	Local methods are typically chosen to compute the Poisson CPD because
	the associated nonconvex optimization problem can be solved efficiently as a 
	sequence of convex subproblems. However, such
	methods are only guaranteed to converge to local optima, which motivates 
	the use of multi-start.%
}
Generalized CP decomposition
(GCP)~\cite{Hong20GeneralizedCanonicalPolyadic,Kolda20StochasticGradientsLargeScale}
incorporates general loss functions into CPD models, including a Poisson 
likelihood-based loss, and stochastic optimization methods. The first stage of 
\textsc{HybridGC} uses GCP with stochastic optimization to form a quick 
approximation which helps the method avoid
local minima that are not the MLE. CP Alternating Poisson Regression
(CPAPR)~\cite{Chi12TensorsSparsityNonnegative} is a deterministic Poisson CPD
method that alternates over a sequence of convex Poisson loss subproblems
iteratively. Previously, in~\cite{Myers21UsingComputationEffectively}, we showed
that CPAPR is performant and can compute accurate approximations to the MLE with
higher probability than GCP. The second stage uses CPAPR to refine the
approximation from GCP to higher accuracy.

\CHANGED{%
	Global methods are typically avoided for CPD due to their high, often
	prohibitive, cost resulting from slow convergence. Nonetheless, they have
	proven 	to be effective for many other global optimization problems. Simulated
	Annealing (SA)~\cite{Kirkpatrick83OptimizationSimulatedAnnealing,Ingber89VeryFastSimulated}
	is one such technique that can handle high-dimensional, nonlinear cost
	functions with arbitrary boundary conditions and constraints, where controlled,
	iterative improvements to the cost function are used in the search for a better model.
	SA effectively leverages stochastic search to avoid local minima and can be
	followed by deterministic search to refine approximate global solutions. Inspired by
	this combined approach, we propose \textsc{HybridGC} for computing the MLE of the Poisson CPD.
}

\subsection{Restarted CPAPR with \textsc{SVDrop}}
\CHANGED{%
	\textsc{HybridGC} is an improvement over existing methods for Poisson CPD 
	but can still converge to local optimizers that are not the MLE, as we illustrate 
	through experimental results presented in~\cref{sec:hybrid-gc}.
	We have identified a specific algebraic property of the approximate solutions
	computed during iterations of
	\textsc{HybridGC} that can lead to such local optimizers.
	Current methods for computing CPD rely on algebraic operations applied to
	tensors that are unfolded into one or more matrix representations.
	In all cases we have explored and are considered in this paper, one or more
	singular values of at least one of these matrix representations of the
	approximate CPD tensor solution computed using \textsc{HybridGC} or 
	one of the local methods drops to nearly zero (i.e., to zero within machine 
	precision) and leads to convergence to one of these local optimizers.
}

\CHANGED{%
	We introduce the parameterized \textsc{SVDrop} heuristic to help identify
	drops of singular values of unfolded CPD solutions during iterations of CPAPR
	(i.e., the iterative refinement method used in the \textsc{HybridGC} method
	described above).  Combining this heuristic with a multi-start strategy, we
	propose the Restarted CPAPR with \textsc{SVDrop} method for Poisson CPD to increase the
	chances of converging to the MLE while reducing the computational cost involved with
	converging to local optima across the multiple initializations.
}

\subsection{Organization}
In \cref{sec:background}, we introduce notation, provide the necessary
background, and discuss related work. In \cref{sec:error}, we formalize several
metrics to compare CPD methods, some of which we used previously
in~\cite{Myers21UsingComputationEffectively}. In \cref{sec:data}, we describe
the data used in numerical experiments. In \cref{sec:hybrid-gc}, we introduce
Hybrid GCP-CPAPR (\textsc{HybridGC}). In \cref{sec:cpapr-svdrop}, we
introduce Restarted CPAPR with \textsc{SVDrop}. In our experiments, we
demonstrate that both methods often improve the likelihood of convergence to the
MLE, thereby reducing excessive computations when compared to multi-start where
the local methods are standalone solvers. In \cref{sec:conc}, we propose future
work.

\section{Background and related work}
\label{sec:background}

\subsection{Notation and conventions}
\label{sec:background:notation}
The set of real numbers and integers are denoted as $\R$ and $\Z$, respectively.
The real numbers and integers restricted to nonnegative values are denoted as
$\R_+$ and $\Z_+$, respectively. The \emph{order} of a tensor is the number of
\textit{dimensions} or \textit{ways}. Each tensor dimension is called a
\textit{mode}. A scalar (tensor of order zero) is represented by a lowercase
letter, e.g., $x$. A bold lowercase letter denotes a vector (tensor of order
one), e.g., $\vec{v}$. A matrix (tensor of order two) is denoted by a bold
capital letter, e.g., $\mat{A} \in \R^{m \times n}$. Tensors of order three and
higher are expressed with a bold capital script letter, e.g., $\tns{X} \in \R^{m
		\times n \times p}$.  Values \emph{computed}, \emph{approximated}, or
\emph{estimated} are typically written with a hat---e.g., $\tns{\widehat{M}} \in
	\R^{m \times n \times p}$.%

The $i$-th entry of a vector $\vec{v}$ is denoted $v_i$, the $(i,j)$ entry of a
matrix $\mat{M}$ is denoted $m_{ij}$, and the $(i,j,k)$ entry of a three-way
tensor $\tns{X}$ is denoted $x_{ijk}$. \emph{Fibers} are the higher-order
analogue of matrix rows and columns. Indices are integer values that range from
1 to a value denoted by the capitalized version of the index variable, e.g., $i
	= 1,\ldots,I$.  We use MATLAB-style notation for subarrays formed from a subset
of indices of a vector, matrix, or tensor mode.  We use the shorthand $i_j
	\colon i_k$ when the subset of indices forming a subarray is the range
$i_j,\ldots,i_k$. The special case of a colon $\colon$ by itself indicates all
elements of a mode, e.g., the $j$-th column or mode-1 fiber of the matrix
$\mat{A}$ is $\mat{A}(\colon,j) = \mat{A}(i_1 \colon i_{I}, j)$. We use the
\textit{multi-index}
\begin{equation}
	\label{eq:cp:multi-index} \mind \coloneqq (i_1, i_2, \ldots, i_d)
	\eqwith i_j \in \{1, 2,\ldots,I_j\} \eqfor j = 1,\ldots,d,
\end{equation}
as a convenient shorthand for the $(i_1, i_2, \ldots, i_d)$ entry of a $d$-way
tensor.

Superscript $T$ denotes non-conjugate matrix transpose. We assume vectors
$\vec{u}$ and $\vec{v}$ are column vectors so that $\vec{u}^T\vec{v}$ is an
inner product of vectors and $\vec{u}\vec{v}^T$ is an outer product of vectors.
We also denote outer products of vectors as $\vec{u} \circ \vec{v} = \vec{u}
	\vec{v}^T$, which is especially useful when describing the $d$-way outer
products of $d$ vectors for $d \geq 2$. The number of matrix or tensor non-zero
elements is denoted $\nnz(\cdot)$; conversely, the number of zeros in a matrix
or tensor is denoted $\nz(\cdot)$.

\subsection{Matricization: transforming a tensor into a matrix}
\textit{Matricization}, as defined
in~\cite{Kolda09TensorDecompositionsApplications}, also known as
\textit{unfolding} or \textit{flattening}, is the process of reordering the
elements of a $d$-way array into a matrix. The mode-$n$ matricization of a
tensor $\tns{X} \in \R^{I_1 \times I_2 \times \dots \times I_d}$, denoted
$\mat{X}_{(n)}$, arranges the mode-$n$ fibers to be the columns of the resulting
matrix.

%

\subsection{Canonical polyadic decomposition}
\label{sec:background:cp}
The canonical polyadic decomposition (CPD) represents a tensor as a finite sum
of rank-one outer products, a generalization of the matrix singular value
decomposition (SVD) to tensors. One major distinction is that there are no
orthogonality constraints on the vectors of the CPD model. Thus we treat the
matrix SVD as a special case of the CPD. Nonetheless, low-rank CP decompositions
are appealing for reasons similar to those of the low-rank SVD, including
dimensionality reduction, compression, de-noising, and more. Interpretability of
CP decompositions on real problems is well-documented, with applications
including exploratory temporal data analysis and link
prediction~\cite{Dunlavy11TemporalLinkPrediction},
chemometrics~\cite{Mocks88TopographicComponentsModel}, neuroscience
\cite{Andersen03PracticalAspectsPARAFAC}, and social network and web link
analysis~\cite{Kolda05HigherorderWebLink,
	Kolda09TensorDecompositionsApplications}.

One particular application of interest is when the tensor data are counts. In
this case, a common modeling choice is to assume that the data follow a Poisson
distribution so that statistical methods, like maximum likelihood estimation,
can be applied to the analysis. One key challenge for computing the Poisson CPD
is that a low-rank CP tensor model of Poisson parameters must satisfy certain
nonnegativity and stochasticity constraints. In the next few sections we cover
the details of the low-rank CP tensor models of Poisson parameters and
decompositions which are the focus of this work.

\subsection{Low-rank CP tensor model}
Assume $\tns{X}$ is a $d$-way tensor of size $I_1 \times \dots \times I_d$. The
tensor $\tns{X}$ is rank-one if it can be expressed as the outer product of $d$
vectors, each corresponding to a mode in $\tns{X}$, i.e.,
\begin{equation}
	\label{eq:cp:rank_one}
	\tns{X} = \vec{a}_1 \circ \vec{a}_2 \circ \dots \circ \vec{a}_d.
\end{equation}
More broadly, the \emph{rank} of a tensor $\tns{X}$ is the smallest number of
rank-one tensors that generate $\tns{X}$ as their
sum~\cite{Kolda09TensorDecompositionsApplications}. We concentrate on the
problem of approximating a tensor of data with a low-rank CP tensor model, i.e.,
the sum of relatively few rank-one tensors.

Let $\bm{\lambda} = [\lambda_1, \lambda_2, \ldots, \lambda_d]\in\R^d$ be a
vector of scalars and let $\mat{A}_1 \in\R^{I_1 \times R}$, $\mat{A}_2 \in
	\R^{I_2\times R}$, $\ldots,$ $\mat{A}_d \in \R^{I_d \times R}$ be matrices. The
\emph{rank-$R$ canonical polyadic (CP) tensor model of
  $\tns{X}$}~\cite{Hitchcock27ExpressionTensorPolyadic} is:
\begin{equation}
	\label{eq:cp}
	\tns{X} \approx \tns{M}
	= \ktns{ \bm{\lambda}; \mat{A}_1, \ldots, \mat{A}_d }
	\coloneqq \sum_{r=1}^R \lambda_r \mat{A}_1(\colon,r) \circ \dots
	\circ \mat{A}_d(\colon,r).
\end{equation}
Each $\mat{A}_k \in \R^{I_k \times R}$ is a \emph{factor matrix} with $I_k$ rows
and $R$ columns. The $j$-th \emph{component} of the mode-$k$ factor matrix is
the column vector $\mat{A}_k(\colon,j)$. We refer to the form $\tns{M} =
	\llbracket \bm{\lambda}; \, \mat{A}_1, \ldots, \mat{A}_d \rrbracket$ as a
\emph{Kruskal tensor}.

\subsection{Computing the Poisson CPD for count data}
\label{sec:background:cp:nntf}
	We focus on an application where all of the entries in a data tensor are counts.
  For the remainder of this work, let $\tns{X} \in \Z_+^{I_1 \times
  \dots \times I_d}$ be a $d$-way tensor of nonnegative integers, let
  $\tns{M}$ be a CP tensor model of the form~\cref{eq:cp}, and assume
  the following about $\tns{X}$:
  \CHANGED{
		\begin{itemize}
			\item each $x_{\mind} \in \tns{X}$ is sampled from a Poisson distribution
				with parameter $m_{\mind}\in \tns{M}$, and
			\item the tensor $\tns{X}$ has low-rank structure (i.e., $R <
				\sqrt[d]{\prod_{k=1}^d I_k}$~\cite{Kolda20StochasticGradientsLargeScale}).
		\end{itemize}%
	}
	Chi and Kolda showed in~\cite{Chi12TensorsSparsityNonnegative} that
	under
	\CHANGED{%
		the assumptions above%
  }
	a \emph{Poisson CP tensor model} is an effective low-rank
	approximation of
	\CHANGED{%
		a tensor of counts, $\tns{X}$, and presented an algorithm for computing
    $\tns{M}$.%
  }
  The Poisson CP tensor model has shown to be
	valuable in analyzing latent patterns and relationships in count data across
	many application areas, including food
	production~\cite{Bro98MultiwayAnalysisFood}, network
	analysis~\cite{Dunlavy11MultilinearAlgebraAnalyzing,
		Baskaran19EnhancingNetworkVisibility}, term-document
	analysis~\cite{Chew07CrosslanguageInformationRetrieval,
		Henretty18TopicModelingAnalysis}, email analysis~\cite{Bader08SurveyTextMining},
	link prediction~\cite{Dunlavy11TemporalLinkPrediction}, geospatial
	analysis~\cite{Ezick19CombiningTensorDecompositions,
		Henretty17QuantitativeQualitativeAnalysis}, web page
	analysis~\cite{Kolda06TOPHITSModelHigherorder}, and phenotyping from electronic
	health records~\cite{Ho14LimestoneHighthroughputCandidate,
		Ho14MarbleHighThroughputPhenotyping,
	Henderson17GraniteDiversifiedSparse}.

	One numerical approach to fit a low-rank Poisson CP tensor model to data is
	\emph{tensor maximum likelihood estimation}, which has proven to be successful in practice.
	Computing the Poisson CPD via tensor maximum likelihood estimation is equivalent to
	solving the following nonlinear, nonconvex optimization problem:
	\begin{equation}
		\label{eq:hybrid-gc:background:nll}
		\min_{\tns{M}} f\left( \tns{X}, \tns{M} \right) =
		\min \sum_{\mind} m_{\mind} - x_{\mind} \log
		m_{\mind},
	\end{equation}
	where $\mind$ is the multi-index~\cref{eq:cp:multi-index}, $x_{\mind} \geq 0$ is
	an entry in $\tns{X}$, and $m_{\mind} > 0$ is a parameter in the Poisson CP
	tensor model $\tns{M}$. The function $f(\tns{X},\tns{M})$
	in~\cref{eq:hybrid-gc:background:nll} is the negative of the log-likelihood of
	the Poisson distribution (omitting the constant
	$\sum_{\mind}\log\left(x_{\mind}!\right)$
	term)~\cite{Rodriguez07PoissonModelsCount}. We will refer to it simply as
\emph{Poisson loss}.

In contrast to linear maximum likelihood
estimation~\cite{Myung03TutorialMaximumLikelihood}, where a single parameter is
estimated using multiple data instances, tensor maximum likelihood estimation
fits a single parameter in an approximate low-rank model to a single data
instance. Within the tensor context, low-rank structure means that multiple
instances in the data are linked to a single model parameter, a type of
multilinear maximum likelihood estimation. This distinction is not made anywhere
else in the literature, to the best of our knowledge.

Much of the research associated with computing the low-rank Poisson CPD via
tensor maximum likelihood estimation has focused on local
methods~\cite{Chi12TensorsSparsityNonnegative,
	Hansen15NewtonbasedOptimizationKullback, Hong20GeneralizedCanonicalPolyadic,
	Kolda20StochasticGradientsLargeScale}, particularly with respect to
computational performance~\cite{Teranishi20SparTenLeveragingKokkos,
	Myers20ParameterSensitivityAnalysis,
	Phipps19SoftwareSparseTensor,Baskaran19FastScalableDistributed,
	Letourneau18ComputationallyEfficientCP, Baskaran17MemoryefficientParallelTensor,
	Baskaran14LowoverheadLoadbalancedScheduling}. Many of the current local methods
for Poisson CPD can be classified as either an
\emph{alternating}~\cite{Carroll70AnalysisIndividualDifferences,
	Harshman70FoundationsPARAFACProcedure} or an
\emph{all-at-once}~\cite{Acar11ScalableOptimizationApproach,
	Acar11AllatonceOptimizationCoupled,Phan13LowComplexityDamped} optimization
method.

Alternating local methods iteratively solve a series of subproblems by fitting
each factor matrix sequentially while the remaining factor matrices are held
fixed. These methods are a form of coordinate descent
(CD)~\cite{Wright15CoordinateDescentAlgorithms}, where each factor matrix is a
block of components that is fit while the remaining component blocks (i.e.,
factor matrices) are left unchanged. Since each block corresponds to a
lower-dimensional problem, alternating tensor methods employ block CD
iteratively to solve a series of easier problems.  CP Alternating Poisson
Regression (CPAPR) was introduced by Chi and Kolda
in~\cite{Chi12TensorsSparsityNonnegative} as a nonlinear Gauss-Seidel approach
to block CD that uses a fixed-point majorization-minimization algorithm called
\emph{Multiplicative Updates (CPAPR-MU)}. At the highest level, the CPAPR
algorithm performs an \textit{outer iteration} where optimizations are applied
on each mode in an alternating fashion. An \emph{inner iteration} is an
optimization using multiplicative updates applied to a subset of variables
corresponding to an individual mode. Inner iterations are performed until the
convergence criterion is satisfied for a mode or up to the maximum allowable
number, $l_{max}$. Outer iterations are performed until the convergence
criterion is satisfied for the whole model or up to the maximum allowable
number, $k_{max}$. The convergence criterion is based on the
\emph{Karush-Kuhn-Tucker (KKT)} conditions, necessary conditions for convergence
to a local minimum in nonlinear optimization. A local minimizer that satisfies
the KKT conditions is called a \emph{KKT point}.

Hansen et al. in~\cite{Hansen15NewtonbasedOptimizationKullback} presented two
Newton-based, active set gradient projection methods using up to second-order
information, \emph{Projected Damped Newton (CPAPR-PDN)} and \emph{Projected
	Quasi-Newton (CPAPR-PQN)}. Moreover, they provided extensions to these methods
where each component block of the CPAPR minimization can be further separated
into independent, highly-parallelizable row-wise subproblems; these methods are
\emph{Projected Damped Newton for the Row subproblem (CPAPR-PDNR)}
and \emph{Projected Quasi-Newton for the Row subproblem (CPAPR-PQNR)}.

One outer iteration of all-at-once optimization methods updates all optimization
variables simultaneously. The Generalized Canonical Polyadic decomposition
algorithm (GCP)~\cite{Hong20GeneralizedCanonicalPolyadic} is a gradient descent
method based on a generic formulation of first derivative information for
arbitrary loss functions to compute the CPD via tensor maximum likelihood
estimation. The original GCP method has two variants: 1) deterministic, which
uses limited-memory quasi-Newton optimization (L-BFGS) and 2) stochastic, which
supports gradient descent (SGD),
AdaGrad~\cite{Duchi11AdaptiveSubgradientMethods}, and
Adam~\cite{Kingma15AdamMethodStochastic} optimizations. The stochastic variants
perform loss function and gradient computations on samples of the input data
tensor so that the search path is computed from estimates of these values. We
focus here on GCP-Adam~\cite{Kolda20StochasticGradientsLargeScale}, which
applies Adam for scalability.

More generally, we focus on the GCP and CPAPR families of tensor maximum
likelihood-based local methods for Poisson CPD for the following reasons:
\begin{enumerate}

	\item \emph{Existing Theory}: Method convergence, computational costs, and
		memory demands are well-understood.
		\CHANGED{%
		  See~\cite[Table
		  1]{Myers20ParameterSensitivityAnalysis} for a summary
		  of convergence results and computational cost,~\cite[\S
		  5.4]{Chi12TensorsSparsityNonnegative} for CPAPR
		  storage results,
		  and~\cite{Kolda20StochasticGradientsLargeScale} for
		  GCP-Adam storage results.
		}

	\item \emph{Available Software}: High-level MATLAB implementing both
	      families is available in Tensor Toolbox for MATLAB
	      (TTB)\footnote{\url{https://gitlab.com/tensors/tensor_toolbox}.}~\cite{Bader17MATLABTensorToolbox,Bader08EfficientMATLABComputations}.
	      A Python version is available in
	      pyttb.\footnote{\url{https://github.com/sandialabs/pyttb}.}
	      High
	      performance C++ code that leverages the Kokkos hardware abstraction
	      library~\cite{Kokkos} to provide parallel computation on diverse
	      computer architectures (e.g., x86-multicore, GPU, etc.) is available
	      with SparTen\footnote{\url{https://github.com/sandialabs/sparten}.}
	      for CPAPR~\cite{Teranishi20SparTenLeveragingKokkos} and
	      Genten\footnote{\url{https://gitlab.com/tensors/genten}.} for
	      GCP~\cite{Phipps19SoftwareSparseTensor}. Additional open-source
	      software for MATLAB includes N-Way
	      Toolbox~\cite{Andersson00NwayToolboxMatlab} and
	      Tensorlab~\cite{Vervliet16TensorlabNumericalOptimization}. Commercial
	      software includes ENSIGN Tensor Toolbox~\cite{Baskaran22ENSIGN}.
\end{enumerate}

\subsection{Error in computing the CPD using multi-start}
\label{sec:error}

Local methods seek local minima. We apply them to global optimization problems
by using a multi-start strategy~\cite{Gyorgy11EfficientMultiStartStrategies,
  Marti18HandbookHeuristics} where a set of approximations are computed from many
random starting points in the feasible domain of the problem. Our methodology is
to generate $N$ random Poisson CP tensor models as initial guesses and compute
$N$ rank-$R$ Poisson CP tensor approximations starting from each initial guess,
which we refer to as \emph{multi-start}. From this set, we choose the ``best''
local minimizer---i.e., the approximation that
minimizes~\cref{eq:hybrid-gc:background:nll}---as the approximation to the
global optimizer.  In turn, the effectiveness of a given method is determined in
part by the probability it will converge to a solution approximating the global
optimizer over all $N$ starting points.

We define several quantities that we will use to compare the effectiveness of a given
method in computing a model that minimizes~\cref{eq:hybrid-gc:background:nll}.
Let $\tns{X}$ be a $d$-way data tensor with dimensions $I_1, \ldots, I_d$. Let
$\mathcal{S} = \{\tns{\widehat{M}}^{(1)},\ldots,\tns{\widehat{M}}^{(N)}\}$ be a
set of rank-$R$ Poisson CP tensor approximations such that $\lvert \mathcal{S}
  \rvert = N$. Let $\tns{M}^*$ denote the \emph{maximum likelihood estimator
  (MLE)}, i.e., the global minimizer of~\cref{eq:hybrid-gc:background:nll}. In
general, \CHANGED{the global minimizer is unknown; however, it has been shown
to exist when $f(\tns{X},\tns{M})$ is finite, which is the case when
$m_{\mind}>0$~\cite{Chi12TensorsSparsityNonnegative}}.
As a result, we aim to recover the \emph{empirical MLE},
$\tns{\widehat{M}}^*$: the rank-$R$ Poisson CP tensor model that is the best
approximation to $\tns{M}^*$. We specify the \emph{empirical MLE restricted to
  $\mathcal{S}$}, i.e., $\tns{\widehat{M}}_{\mathcal{S}}^* \equiv
  \tns{\widehat{M}}^* \in \mathcal{S}$, as the best approximation to the MLE from
$\mathcal{S}$:
\begin{equation}
  \tns{\widehat{M}}_{\mathcal{S}}^* = \{ \tns{\widehat{M}}^{(j)} \in
  \mathcal{S} \mid
  f(\tns{X},\tns{\widehat{M}}^{(j)}) \leq f(\tns{X},
  \tns{\widehat{M}}^{(k)} ), \, k = 1, \ldots, |\mathcal{S}|, \, j < k \}.
\end{equation}
The condition that $j<k$ guarantees that the set is nonempty in the case of a
tie. We write $\mathcal{S}_A$ when every element in $\mathcal{S}$ was computed
by algorithm $A$. This notation will be useful later on when analyzing results
from different algorithms.

\subsection{An error estimator on the loss function}
The probability that algorithm $A$ converges from any starting point in the
feasible region of~\cref{eq:hybrid-gc:background:nll} to some
$\tns{\widehat{M}}^{(n)} \in \mathcal{S}_A$ such that
$f(\tns{X},\tns{\widehat{M}}^{(n)})$ is within a ball of radius $\epsilon>0$ of
$f(\tns{X}, \tns{\widehat{M}}^*)$ is defined as
\CHANGED{%
  \begin{equation}
    \label{eq:hybrid-gc:error:nll:prob_mle}
    P_A\left( \rho_n
    < \epsilon\right) \quad \text{as} \quad
    n = 1,\ldots,|\mathcal{S}_A|,
  \end{equation}
  where
  \begin{equation}
      \rho_n \coloneqq \frac{ \lvert f(\tns{X},
        \tns{\widehat{M}}^{(n)}) - f(\tns{X},\tns{\widehat{M}}^*) \rvert }{ \lvert
        f (\tns{X},\tns{\widehat{M}}^*) \rvert }
      \end{equation}
  is the relative error in Poisson loss.
  We can only estimate $P_A$ since $\mathcal{S}_A$ is finite in practice.
  A computable estimator to~\cref{eq:hybrid-gc:error:nll:prob_mle} is%
}
\begin{equation}
  \label{eq:hybrid-gc:error:nll:prob_mle_est}
    \widehat{P} \left( \mathcal{S}_A, \epsilon \right) = \frac{\# \tns{\widehat{M}}^{(n)} \in \mathcal{S}_A
      \text{ for which } \rho_n <  \epsilon }{\lvert \mathcal{S}_A \rvert}.
\end{equation}
\Cref{eq:hybrid-gc:error:nll:prob_mle_est} is a conservative estimator since it
does not account for solutions which may be closer to $\tns{M}^*$ but are more
than $\epsilon$-distance from $\tns{\widehat{M}}^*$. We omit $\mathcal{S}_A$ and
write only $\widehat{P}(\epsilon)$ when the method is clear from the context.

\subsection{An error estimator on the algebraic structures}
We define two measures of approximation error quantifying the structural
similarity between two Kruskal tensors based on their algebraic properties
called \emph{factor match score}~\cite{Korth75DistributionChanceCongruence,
  Korth76ProcrustesMatchingCongruence, Chi12TensorsSparsityNonnegative,
  Kolda20StochasticGradientsLargeScale}.

\subsubsection{Factor match score (FMS)}

FMS is the maximum sum of cosine similarities over all permutations of the
column vectors of all the factor matrices between two Kruskal tensors,
$\tns{M}_1 = \llbracket \mat{\lambda}^{\mat{A}}; \, \mat{A}_1, \ldots, \mat{A}_d
\rrbracket$ and $\tns{M}_2 = \llbracket \mat{\lambda}^{\mat{B}}; \, \mat{B}_1,
\ldots, \mat{B}_d \rrbracket$:
\CHANGED{
  \begin{equation}
    \label{eq:hybrid-gc:error:fms:def_fms}
    \begin{gathered}
      \FMS(\tns{M}_1,\tns{M}_2) =
      \max_{\pi( \cdot )}
      \frac{1}{R} \sum_{r=1}^R
      \left( 1 -
      \frac{\abs{ \xi_r - \zeta_r }}
      {\max \{ \xi_r, \zeta_r \}}
      \right)
      \prod_{n=1}^d
      \frac{
        \mat{A}_n (\colon, j)^T
        \mat{B}_n (\colon,\pi(j))
      }
      {
        \norm{\mat{A}_n (\colon, j)}
        \norm{\mat{B}_n (\colon, \pi(j))}
      },
      \\
      \eqwhere
      \xi_r = \lambda_r^{\mat{A}} \prod_{n=1}^d
      \norm{ \mat{A}_n(\colon,r) } \quad \text{and} \quad
      \zeta_r = \lambda_r^{\mat{B}} \prod_{n=1}^d
      \norm{ \mat{B}_n(\colon,r) }.
    \end{gathered}
  \end{equation}
  The permutation $\pi(\cdot)$%
}
reorders the columns of the
factor matrices of $\tns{M}_2$ to maximize the number of columns that are
correctly identified.
An FMS of $1$ indicates collinearity among the columns of all factor matrices
and thus a perfect match between the two Kruskal tensors. As
in~\cite{Lorenzo-Seva06TuckerCongruenceCoefficient}, we say $\tns{M}_1$ and
$\tns{M}_2$ are \emph{similar} if $\FMS(\tns{M}_1,\tns{M}_2) \geq 0.85$ and
\emph{equal} if $\FMS(\tns{M}_1,\tns{M}_2) \geq 0.95$, which are common values
used to define acceptable matches in recent
work~\cite{Chi12TensorsSparsityNonnegative,
  Hansen15NewtonbasedOptimizationKullback, Kolda20StochasticGradientsLargeScale}.
FMS is a particularly useful measure of the effectiveness of a method in
relating the low-rank structure of an approximation to that of a known model.
Using FMS, we estimate the probability that a method computes models with the
same algebraic structure as the empirical MLE. We formalize this now.
\subsubsection{Probability of similarity}
For each computed solution $\tns{\widehat{M}}^{(n)} \in \mathcal{S}$,
$n=1,\ldots,|\mathcal{S}|$, define an indicator function
$\psi_n(\tns{M},\tns{\widehat{M}}^{(n)},t)$ that is 1 when the $n$-th model has
$\FMS (\tns{M},\tns{\widehat{M}}^{(n)}) \geq t$ and 0 otherwise; i.e.,
\begin{equation}
  \label{eq:hybrid-gc:error:fms:indicator}
  \psi_n(\tns{M},\tns{\widehat{M}}^{(n)},t) =
  \begin{cases*}
    1, \quad \text{if } \FMS(\tns{M}, \tns{\widehat{M}}^{(n)}) \geq t \\
    0, \quad \text{otherwise}.
  \end{cases*}
\end{equation}
We use \cref{eq:hybrid-gc:error:fms:indicator} in our discussions below to
quantify the \emph{fraction over $N$ solves with FMS greater than $t$},
\begin{equation}
  \label{eq:hybrid-gc:error:fms:prob_mle_est}
  \begin{gathered}
    \Psi(\tns{M},\mathcal{S},t) =
    \frac{1}{N}\sum_{n=1}^{|\mathcal{S}|} \psi_n(\tns{M}, \tns{\widehat{M}}^{(n)},t), \\
    \eqwhere \tns{\widehat{M}}^{(n)} \in \mathcal{S}, \, \forall n \in \{1, \ldots, |\mathcal{S}|\}.
  \end{gathered}
\end{equation}

\subsection{Related work}
There are other approaches in the literature that seek to fit models with other
distributions in the exponential family or that use other algorithms to estimate
parameters. Alternating least squares methods are relatively easy to implement
and effective when used with LASSO-type
regularization~\cite{Friedlander08ComputingNonnegativeTensor,
	Bazerque13InferencePoissonCount}. The method of Ranadive et
al.~\cite{Ranadive21AllOnceCP}, CP-POPT-DGN, is an all-at-once active set
trust-region gradient-projection method. CP-POPT-DGN is functionally very
similar to CPAPR-PDN. Whereas CP-POPT-DGN computes the search direction via
preconditioned conjugate gradient (PCG), CPAPR-PDNR computes the search
direction via Cholesky factorization. The most significant differences are: 1)
CP-POPT-DGN is all-at-once whereas all CPAPR methods are alternating and 2)
CPAPR can take advantage of the separable row subproblem formulation to achieve
more fine-grained parallelism. The Generalized Gauss-Newton method of
Vandecapelle et al.~\cite{Vandecappelle20SecondorderMethodFitting} follows the
GCP framework to fit arbitrary non-least squares loss via an all-at-once
optimization and trust-region-based Gauss-Newton approach.  Hu et
al.~\cite{Hu15ScalableBayesianNonnegative, Hu15ZeroTruncatedPoissonTensor}
re-parameterized the Poisson regression problem to leverage Gibbs sampling and
variational Bayesian inference to account for the inability of CPAPR to handle
missing data. Other problem transformations include probabilistic likelihood
extensions via Expectation Maximization~\cite{Rai15ScalableProbabilisticTensor,
	Huang17KullbackLeiblerPrincipalComponent} and a Legendre
decomposition~\cite{Sugiyama18LegendreDecompositionTensors} instead of a CP
decomposition.

\CHANGED{%
	\Cref{sec:cpapr-svdrop} presents the \textsc{SVDrop} heuristic that can identify
	approximate solutions that will converge to non-MLE solutions (i.e., local optimizers)
	for Poisson CPD early in the iterative process. This heuristic is based on spectral
	properties of unfoldings of the approximate solution, and previous work has also
	attempted to characterize non-optimal solutions for the CPD problem based on
	related properties.
	Early work by Kruskal et al.~\cite{Kruskal89How3MFAData} studied two-factor
	degeneracies (2FD), where two components of a CPD model are highly correlated in
	all three modes. Mitchell and Burdick introduced a test for 2FD
	in~\cite{Mitchell94SlowlyConvergingParafac}, which  uses an early definition of
	FMS to identify 2FD between successive iterates. More recent
	work~\cite{Giordani16RemediesDegeneracyCandecomp} adds constraints (e.g.,
	orthogonality constraints, ridge regression, SVD penalty) to the CPD problem to
	avoid 2FD.
}

\CHANGED{%
	Breiding and Vannieuwenhoven proposed a formulation of CPD
	as a Riemannian optimization problem, described spectral properties
	of approximate solutions to this reformulation that lead to ill-conditioning, and
	proposed a restart technique to escape from regions of ill-conditioning where
	iterations can stagnate~\cite{BrVa18siopt,BrVa18simax,BrVa18amletters}.
	Their work most closely resembles the ideas presented here,
	including identifying spectral properties of iterations that lead to non-optimal solutions
	and restarting sub-optimal iterations before full convergence is attained to reduce
	computational cost. However, their approach was developed for CPD with Gaussian data
	assumptions, including the Riemannian optimization objective function and derivatives
	(whose spectral properties are used for identifying ill-conditioning), a Riemannian
	optimization retraction function based on a low-rank Tucker decomposition,
	ST-HOSVD~\cite{VaVaMe12} that has not been extended directly to Poisson CPD, and a
	restart method based on local perturbations of iterates in regions of ill-conditioning. Our
	work here differs in that we focus exclusively on methods applicable to Poisson CPD, leverage
	spectral properties of the unfoldings of iterates (rather than of the Hessian of the objective
	function at iterates, thus not requiring second-order derivative information that may not be
	readily available), and use a restart technique that does not rely on the current iterates
	(as empirical evidence with Poisson CPD solvers indicated that such local perturbations
	did not consistently lead to convergence to MLE solutions). Moreover, the results presented here are
	empirical, whereas Breiding and Vannieuwenhoven focused on both theoretical analysis and empirical
	evidence to identify ill-conditioning of CPD and propose remedies. Future work could consider extending
	their ideas to better understand the challenges associated with Poisson CPD that we identify and address
	in this paper.%
}

\section{Data examples}
\label{sec:data}

\subsection{\texttt{LowRankSmall}} A synthetic count tensor with
dimensions $4 \times 6 \times 8$, generated rank $r=3$, and 17 nonzero
entries (8.85\% dense).
The sparse tensor is fully provided in~\Cref{sec:data:low-rank-small}
as~\Cref{data:hybrid-gc:low-rank-small-sptensor}. All experiments
featuring this dataset were conducted with Tensor Toolbox for MATLAB
v3.3. Additional implementation details are specified
in~\cref{sec:hybrid-gc:experiments}
and~\cref{sec:cpapr-svdrop:experiments}.

\subsection{\texttt{MedRankLarge}} A synthetic count tensor
 with dimensions $1000 \times 1000 \times 1000$ , generated rank
 $r=20$, and $98026$ nonzero entries (0.009\% dense).
 All experiments featuring this dataset were conducted with SparTen for
 CPAPR and Genten for GCP-Adam on a dual-socket Intel Xeon Gold
 processor with 18 cores per socket. Both tools were compiled for single-node
 parallelism with GCC v8.3.1 and were run using all available OpenMP
 threads. Additional implementation details are specified
 in~\cref{sec:hybrid-gc:experiments}.

\section{Hybrid GCP-CPAPR}
\label{sec:hybrid-gc}

We present Hybrid GCP-CPAPR (\textsc{HybridGC}), an algorithm for Poisson CPD
that estimates the solution of a nonlinear, nonconvex optimization problem by
\emph{approximating a global optimization algorithm through the composition of
local methods}. \textsc{HybridGC} first uses stochasticity to compute a
coarse-grained estimate of the model and then refines the model with a
deterministic method. Our numerical experiments demonstrate the synergy of this
hybrid approach: \textsc{HybridGC} yields an effective algorithm that computes
an approximation to the MLE for Poisson CPD with higher accuracy than the
methods it leverages. The stochastic stage makes \textsc{HybridGC} scalable to
very large problems and the deterministic stage allows the method to exploit
convergence results in~\cite{Chi12TensorsSparsityNonnegative}. To the best of
our knowledge, this is the first work that extends similar approaches in the
matrix case (for example,~\cite{Halko11FindingStructureRandomness}) to low-rank
tensor decompositions in this way.

\subsection{\textsc{HybridGC} method}
Like SA, \textsc{HybridGC} leverages both stochastic and deterministic
optimizations to start from an initial guess, iterate according to
a schedule, and converge to a solution approximating the global
optimizer. Specifically, \textsc{HybridGC} iterates from an initial
guess $\tns{M}_{0}$ via a two-stage optimization between stochastic
and deterministic search to return a Poisson CP tensor model
$\tns{\widehat{M}}$ that is an estimate to $\tns{M}^*$. In the first
stage, the stochastic search method starts from $\tns{M}_{0}$ and
iterates for $j$ iterations to return an intermediate solution,
$\tns{M}_1$. In contrast to SA, which uses random perturbation for the
``heating'' step, we use
GCP-Adam~\cite{Kolda20StochasticGradientsLargeScale} for structured
stochastic optimization. In the second stage, deterministic search
refines $\tns{M}_1$ for $k$ iterations to return $\tns{M}_2$. We use
CPAPR with Multiplicative
Updates~\cite{Chi12TensorsSparsityNonnegative} as our ``cooling''
step. \textsc{HybridGC} returns $\tns{\widehat{M}} = \tns{M}_2$ as an
estimate to the global optimizer, $\tns{M}^*$. The details of
\textsc{HybridGC} are given below in \Cref{alg:hybrid-gc}.%

Presently, our analogue of the temperature schedule is the combination
of the number of GCP epochs and CPAPR outer iterations.
\CHANGED{%
		A GCP \emph{epoch} is comprised of one or more \emph{iterations}. The Adam step parameter and estimate of the objective function value both are updated once per epoch. The stochastic gradient is computed and a descent step is taken in each iteration. A CPAPR \emph{outer iteration} corresponds to one pass over all of the tensor modes. For each mode, the tensor is matricized and one or more \emph{inner iterations} are taken to optimize the factor matrix corresponding to that mode.
	}

In further contrast to SA, we do not include a notion of
acceptance-rejection with respect to newly obtained states; we leave
this to future work. We only consider stochastic search followed by
deterministic search and not the opposite. This is because stochastic
search directions are found using estimates of the objective function
from sample points. Thus it would be possible for the algorithm to
converge to a minimum yet remain marked as \emph{not converged} if the
objective function value were only coarsely estimated. Thus it is
likely that stochastic search would move away from the optimum.

\begin{algorithm}
    \caption{Hybrid GCP-CPAPR}
    \label{alg:hybrid-gc}
    \begin{algorithmic}
        \Function{HybridGC}{tensor $\tns{X}$, rank $r$, initial guess
            $\tns{M}_0$} \State $\tns{M}_1 \gets \textsc{GCP}(\tns{X}, \,r,
            \,\tns{M}_0)$ \State $\tns{M}_2 \gets \textsc{CPAPR}(\tns{X}, \,r,
            \,\tns{M}_1)$ \State \Return model tensor $\tns{\widehat{M}} =
            \tns{M}_2$ as estimate to $\tns{M}^*$ \EndFunction
    \end{algorithmic}
\end{algorithm}

\subsection{Numerical experiments}
\label{sec:hybrid-gc:experiments}

This section presents experiments that were designed to evaluate the algorithm
effectiveness of \textsc{HybridGC} compared to GCP-Adam and CPAPR-MU
as standalone solvers.
\CHANGED{%
	We demonstrate this by performing many independent trials on each of
	two synthetic low-rank Poisson multilinear datasets. For each trial
	on a dataset, we compute three rank-$r$ Poisson CPD approximations with
	GCP-Adam, CPAPR-MU, and HybridGC; the specific parameterizations are
	detailed below in~\cref{sec:hybrid-gc:experiments:proc}.
}
For the remainder of this
work, we refer to GCP and CPAPR without further specifying the optimization
routine. We use $\mathcal{S}_G$, $\mathcal{S}_C$, and $\mathcal{S}_H$ to refer
to the sets of approximations computed with GCP, CPAPR, and HybridGC,
respectively.

We treat \textsc{HybridGC} decompositions as those computed with $j\geq 0$ GCP epochs followed by $k \geq 0$ CPAPR outer iterations.
\CHANGED{%
By default, which we use in our experiments, GCP sets $1 \text{ epoch} = 1000 \text{ iterations}$ and CPAPR sets a maximum of 10 inner iterations per mode per outer iteration.
}

Using the notation set above, we denote the solutions that were
computed with GCP, CPAPR, and \textsc{HybridGC} as
\begin{align*}
	\mathcal{S}_G & = \{ \tns{\widehat{M}}_j \mid
	\tns{\widehat{M}}_j \text{ computed by GCP}\};                \\
	\mathcal{S}_C & = \{ \tns{\widehat{M}}_j \mid
	\tns{\widehat{M}}_j \text{ computed by CPAPR}\}; \text{ and,} \\
	\mathcal{S}_H & = \{ \tns{\widehat{M}}_j \mid
	\tns{\widehat{M}}_j \text{ computed by \textsc{HybridGC}}\}.
\end{align*}

\subsubsection{Procedure}
\label{sec:hybrid-gc:experiments:proc}
Fix an input tensor $\tns{X} \in \R^{I_1 \times \dots \times I_d}$ and a rank
$R$. Specify a maximum work budget $W=J_{max} + K_{max}$ for all methods where
$J_{max} \geq 0$ and $K_{max} \geq 0$ are the maximum allowable number of outer
iterations for GCP and CPAPR, respectively. If $J_{max} = 0$ and $K_{max} = W$,
then \textsc{HybridGC} is equivalent to CPAPR. Conversely, if $J_{max} = W$ and
$K_{max} = 0$, then \textsc{HybridGC} is equivalent to GCP. In this way,
\textsc{HybridGC} generalizes both methods by ``interpolating'' GCP and CPAPR
when both $J_{max} >0$ and $K_{max} > 0$.

Starting from the same random initial guess, we computed $j \in (0,\ldots,W)$
rank-$R$ decompositions with GCP iterating for at most $j$ epochs. GCP
is considered converged when the stochastic gradient learning rate $\alpha$ is
decayed from $10^{-3}$ (the default) to $10^{-15}$. Next, starting from each of
the $W+1$ epoch iterates computed with GCP, we computed $k \in (W, W-1, \ldots, 0)$
rank-$R$ decompositions with CPAPR iterating for at most $k$ outer iterations.
CPAPR is considered converged when the KKT-based criterion is less than or equal
to $10^{-15}$. Since each \textsc{HybridGC} trial produced $W+1$ decompositions,
only the empirical MLE restricted to that trial was chosen for comparison.

\CHANGED{%
	In total, $N=\text{110,266}$ rank $R=3$ decompositions of \texttt{LowRankSmall} were computed with each
	algorithm. Two minimizers were found in this set of solutions and the absolute minimum was chosen
	as the empirical MLE, with $f(\tns{X},\tns{M}) = -26.214062308801800$. Among the $N=100$ rank
	$R=20$ decompositions of \texttt{MedRankLarge} computed with each algorithm, many minimizers were
	found. The decomposition with the absolute minimum over all trials was chosen as the empirical
	MLE. In both cases, we refer to the set of solutions within a ball of radius $\epsilon$ around the
	empirical MLE as having converged to the MLE for some $\epsilon > 0$. In the case of
	\texttt{LowRankSmall}, since we have only one other minimizer, in our discussions we refer
	explicitly to solutions within a ball of $\epsilon$ around the second local minimizer, with
	$f(\tns{X},\tns{M})=-5.838518788084730$, as having converged to \emph{the local minimum} for some
	$\epsilon>0$. We specify $\epsilon$ where it is germane to the discussion.
}

\subsubsection{Comparison on the loss function}
\label{sec:hybrid-gc:experiments:nll}
We now compare the effectiveness of \textsc{HybridGC} as an algorithm for
solving a nonlinear, nonconvex optimization problem by considering the Poisson
loss~\cref{eq:hybrid-gc:background:nll} and the MLE-probability estimator based
on it~\cref{eq:hybrid-gc:error:nll:prob_mle_est}. The empirical MLE is denoted
$\tns{M}_{\mathcal{S}}^*$, with $\mathcal{S} = \mathcal{S}_G \cup \mathcal{S}_C
	\cup \mathcal{S}_H$.

\Cref{tab:hybrid-gc:experiments:nll:prob_mle} presents average behavior about
algorithm effectiveness as estimates of the probability that each method
computed the empirical MLE with relative
error~\cref{eq:hybrid-gc:error:nll:prob_mle_est} less than $\epsilon$ for
\texttt{LowRankSmall} and \texttt{MedRankLarge}. Viewing the average behavior,
we note the following:
\begin{itemize}
	\item For the small dataset with low rank (\texttt{LowRankSmall}),
	      \textsc{HybridGC} and CPAPR had comparable precision at all levels of
	      accuracy. Additional work may be conducted to determine if the
	      differences are statistically significant.
	\item \CHANGED{%
					\textsc{HybridGC} always converged to the MLE when CPAPR
					did. In our experiments with \texttt{LowRankSmall},
					\textsc{HybridGC} converged to the MLE in the same trials as
					CPAPR did plus in an additional 0.42\% of trials.
				}
	\item \CHANGED{%
					\textsc{HybridGC} always converged to the MLE at the same or
					a higher rate than GCP.
					Although there were instances on \texttt{LowRankSmall} where GCP
					converged to the MLE and \textsc{HybridGC} did not,
					the opposite occurred in 13.7\% more trials.
				}
	\item \CHANGED{%
					On both datasets \textsc{HybridGC} had a
					higher probability of getting close to the empirical MLE at high accuracy.
				}

	\item Even for small input tensors with low rank, GCP was virtually
	      incapable of resolving the MLE beyond a coarse-grain approximation.
	      The situation was worse for larger tensors with more components.
\end{itemize}

\begin{table}[!ht]
	\centering
	\caption{Estimate of probability each method computes a solution within
		$\epsilon$-radius of approximate global optimizer.}
	\label{tab:hybrid-gc:experiments:nll:prob_mle}
	\begin{subtable}{0.45\textwidth}
		\centering
		\caption{\texttt{LowRankSmall} dataset ($>110$K trials).}
		\label{tab:hybrid-gc:experiments:nll:prob_mle:data1}
		\begin{tabular}{crrr}
			\toprule
			$\epsilon$ & CPAPR & GCP   & \textsc{HybridGC} \\
			\midrule
			$10^{-1}$  & 0.963 & 0.963 & \textbf{0.967}    \\
			$10^{-2}$  & 0.963 & 0.963 & \textbf{0.967}    \\
			$10^{-3}$  & 0.963 & 0.879 & \textbf{0.967}    \\
			$10^{-4}$  & 0.963 & 0.003 & \textbf{0.967}    \\
			\bottomrule
		\end{tabular}
	\end{subtable}\qquad
	\begin{subtable}{0.45\textwidth}
		\centering
		\caption{\texttt{MedRankLarge} dataset (100 trials).}
		\label{tab:hybrid-gc:experiments:nll:prob_mle:data2}
		\begin{tabular}{crrr}
			\toprule
			$\epsilon$ & CPAPR         & GCP  & \textsc{HybridGC} \\
			\midrule
			$10^{-1}$  & 1.00          & 1.00 & 1.00              \\
			$10^{-2}$  & \textbf{0.46} & 0.04 & \textbf{0.46}     \\
			$10^{-3}$  & 0.03          & 0.00 & \textbf{0.17}     \\
			$10^{-4}$  & 0.00          & 0.00 & \textbf{0.01}     \\
			\bottomrule
		\end{tabular}
	\end{subtable}
\end{table}

\begin{figure}[!ht]
	\centering
	\begin{subfigure}[t]{1.00\textwidth}
		\includegraphics[width=\textwidth]{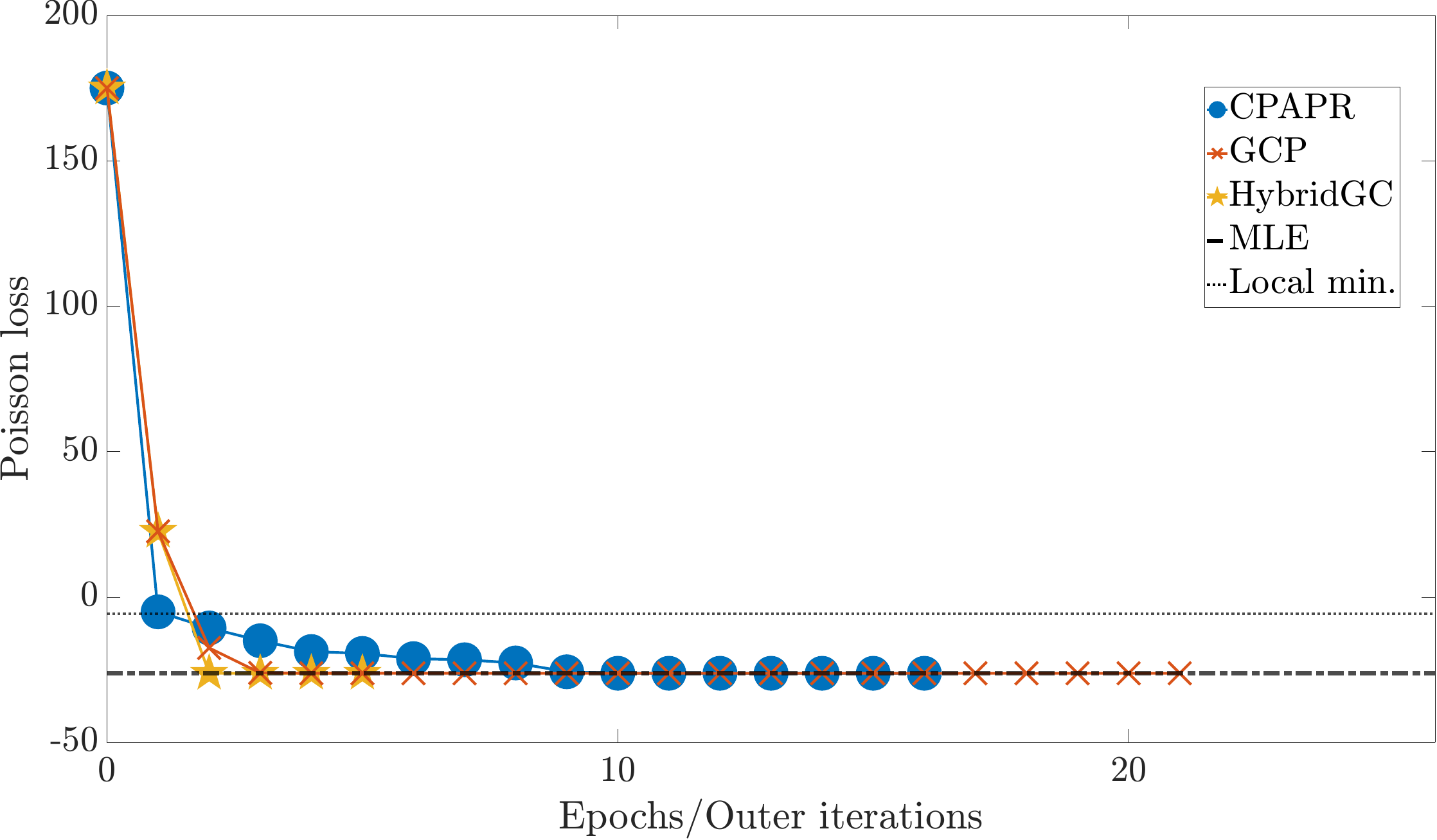}
		\caption{Trial where all methods converged to the MLE.}
		\label{fig:hybrid-gc:experiments:nll:data1:traces:mle}
	\end{subfigure}\\
	\begin{subfigure}[t]{1.00\textwidth}
		\includegraphics[width=\textwidth]{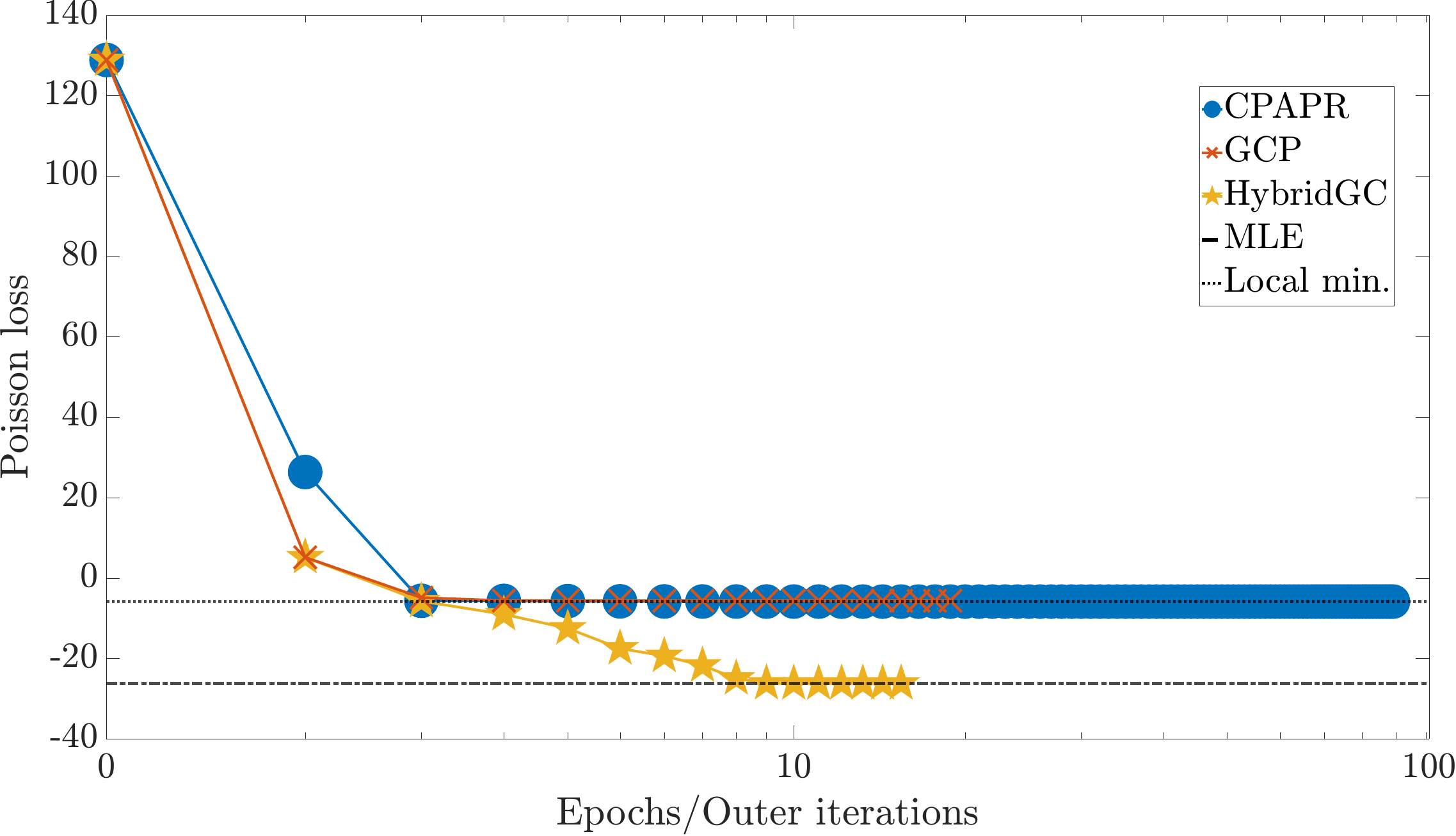}
		\caption{Trial where only \textsc{HybridGC} converged to the MLE.}
		\label{fig:hybrid-gc:experiments:nll:data1:traces:loc}
	\end{subfigure}
	\caption{Examples of two types of behaviors of traces of loss function
		values for GCP, CPAPR, and \textsc{HybridGC} on \texttt{LowRankSmall}.
		\CHANGED{%
			The MLE and the second
			local minimum are shown in both plots for direct comparison.%
		}
	}
	\label{fig:hybrid-gc:experiments:nll:data1:traces}
\end{figure}

\Cref{fig:hybrid-gc:experiments:nll:data1:traces} highlights two
behaviors observed in two different trials among the $N=\text{110,266}$ trials
we ran on \texttt{LowRankSmall}.
\Cref{fig:hybrid-gc:experiments:nll:data1:traces:mle} presents the traces in
Poisson loss function value from one trial where all methods
computed approximations close to the empirical MLE when started from the same
initial guess for \texttt{LowRankSmall}.
\Cref{fig:hybrid-gc:experiments:nll:data1:traces:loc} presents similar traces
except when only \textsc{HybridGC} converged to the MLE and the standalone
solvers, GCP and CPAPR, converged to a different local minimum.
Of all $N=\text{110,266}$ trials, the empirical MLE was computed by \textsc{HybridGC}.
\CHANGED{%
	The MLE and the second local minimum are shown in
	both~\Cref{fig:hybrid-gc:experiments:nll:data1:traces:mle}
	and~\Cref{fig:hybrid-gc:experiments:nll:data1:traces:loc} for direct comparison with the same
	axes.%
}

\begin{figure}[!ht]
	\centering
	\includegraphics[width=\textwidth]{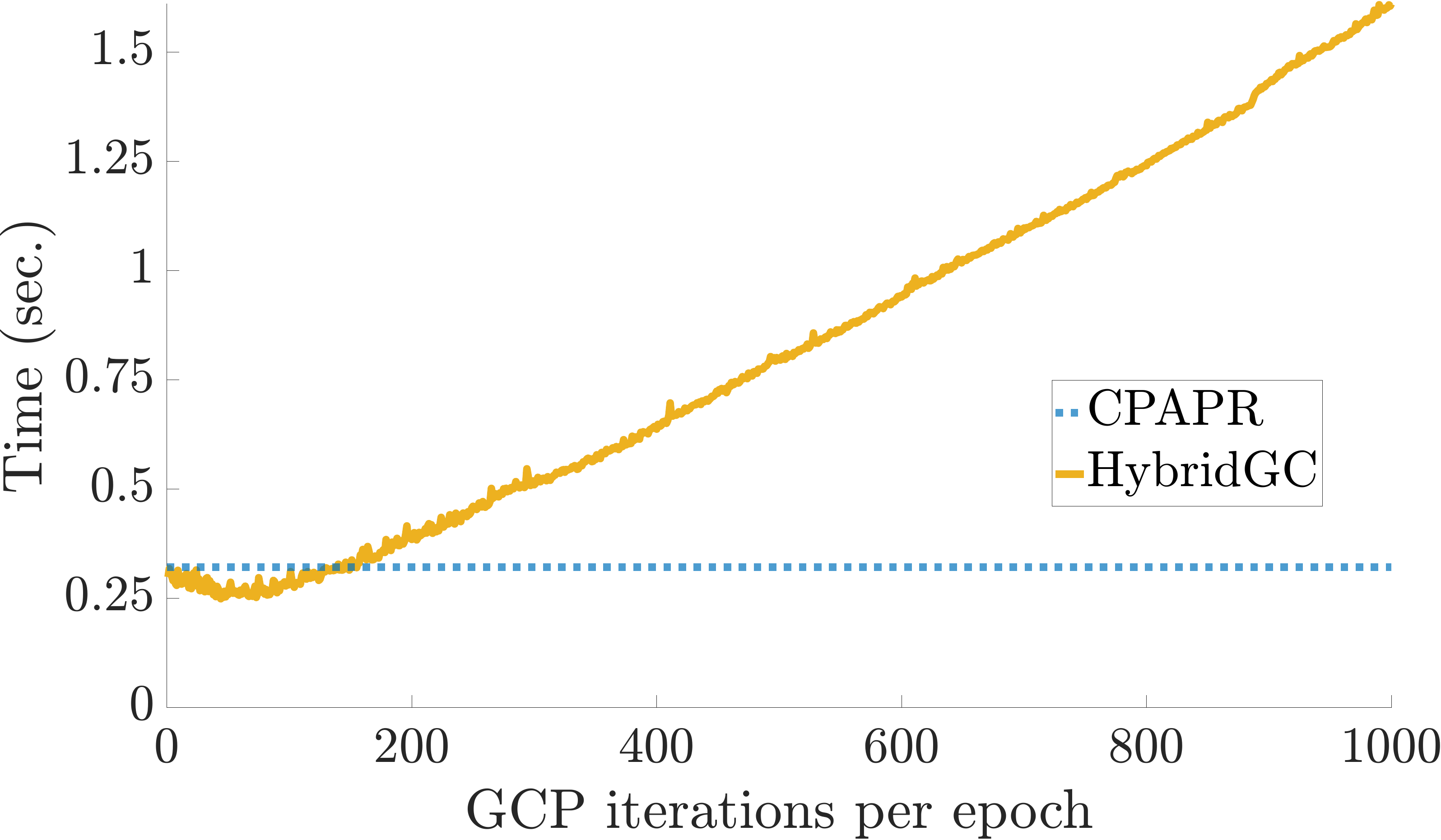}
	\caption{
		\CHANGED{%
			Performance of \textsc{HybridGC} versus CPAPR as a standalone solver where both methods converge to the MLE. The $x$-axis indicates how many iterations were run in the GCP stage before starting the CPAPR stage of \textsc{HybridGC}.%
		}
	}
	\label{fig:hybrid-gc:experiments:nll:data1:timings}
\end{figure}

\CHANGED{%
	\Cref{tab:hybrid-gc:experiments:nll:data1:times} shows the wall clock
	times of the two trials presented in~\Cref{fig:hybrid-gc:experiments:nll:data1:traces}.
	It is clear that GCP performs worse than CPAPR and
	\textsc{HybridGC} for this parameterization regardless of the
	minimizer it converges to.  When CPAPR converged to the MLE, it converged approximately
	roughly $4\times$ faster than \textsc{HybridGC}. Closer
	examination reveals that more than 97\% of \textsc{HybridGC} runtime was spent in
	the GCP stage. We found the default parameterization (i.e., $1 \text{ epoch} = 1000 \text{ iterations}$) to determine the poor performance of \textsc{HybridGC} compared to CPAPR. To see this, we studied the convergence behavior and performance of running the CPAPR stage of \textsc{HybridGC} starting from the approximation computed after every iteration in the first epoch of the GCP stage. Not only did \textsc{HybridGC} converge to the MLE from every starting point but \textsc{HybridGC} outperformed CPAPR in 14.5\% of starts, shown in~\cref{fig:hybrid-gc:experiments:nll:data1:timings}. In the best case, \textsc{HybridGC} was 29\% faster than CPAPR. Since it was not possible to conduct a similar experiment for all trials due to storage and computational limitations, we acknowledge that this result may not hold in general. Nonetheless, this anecdotal evidence suggests, when considered with the overall convergence results described above, that \textsc{HybridGC} with very small epoch sizes in the GCP stage may be highly effective in computing the MLE efficiently with high probability.
}

\begin{table}[t]
	\centering
		\caption{Total time (sec.) for each trial in~\Cref{fig:hybrid-gc:experiments:nll:data1:traces}.}
		\label{tab:hybrid-gc:experiments:nll:data1:times}
		\begin{tabular}{lrrr}
		\toprule
		Trial & CPAPR & GCP & \textsc{HybridGC} \\
		\midrule
		\Cref{fig:hybrid-gc:experiments:nll:data1:traces:mle} & 0.35 & 26.58 & 1.39 \\
		\Cref{fig:hybrid-gc:experiments:nll:data1:traces:loc} & 2.34 & 21.97 & 1.51 \\
		\bottomrule
	\end{tabular}
\end{table}

A few additional remarks:
\begin{itemize}
	\item \CHANGED{%
					In our experiments we observe%
				}
				the benefit of performing some amount of stochastic search
				followed by deterministic search. Owing to performing one
				epoch of stochastic search, \textsc{HybridGC}
				quickly identified the basin of attraction of the MLE and
				converged before GCP and CPAPR in both trials, even in the
				case where the standalone method converges to a local
				minimizer that is different from the MLE. The iteration
				histories of GCP and CPAPR are consistent with results from
				prior work on small
				tensors~\cite{Myers21UsingComputationEffectively}.
	\item The shared behavior among all methods of making only incremental
	      progress before finally converging is a feature of the theoretical
	      convergence properties of each method. Mathematically, CPAPR (in the
	      case of MU) and the deterministic stage of \textsc{HybridGC} converge
	      sublinearly in the basin of attraction to the MLE; GCP (in the case of
	      Adam) converges only linearly at best.
				\CHANGED{%
					See~\cite[Table 1]{Myers21UsingComputationEffectively}
					for details.%
				}
\end{itemize}

\subsubsection{Comparison as algebraic structures}
\label{sec:hybrid-gc:experiments:fms}
\begin{figure}[!ht]
	\centering
	\includegraphics[width=\textwidth]{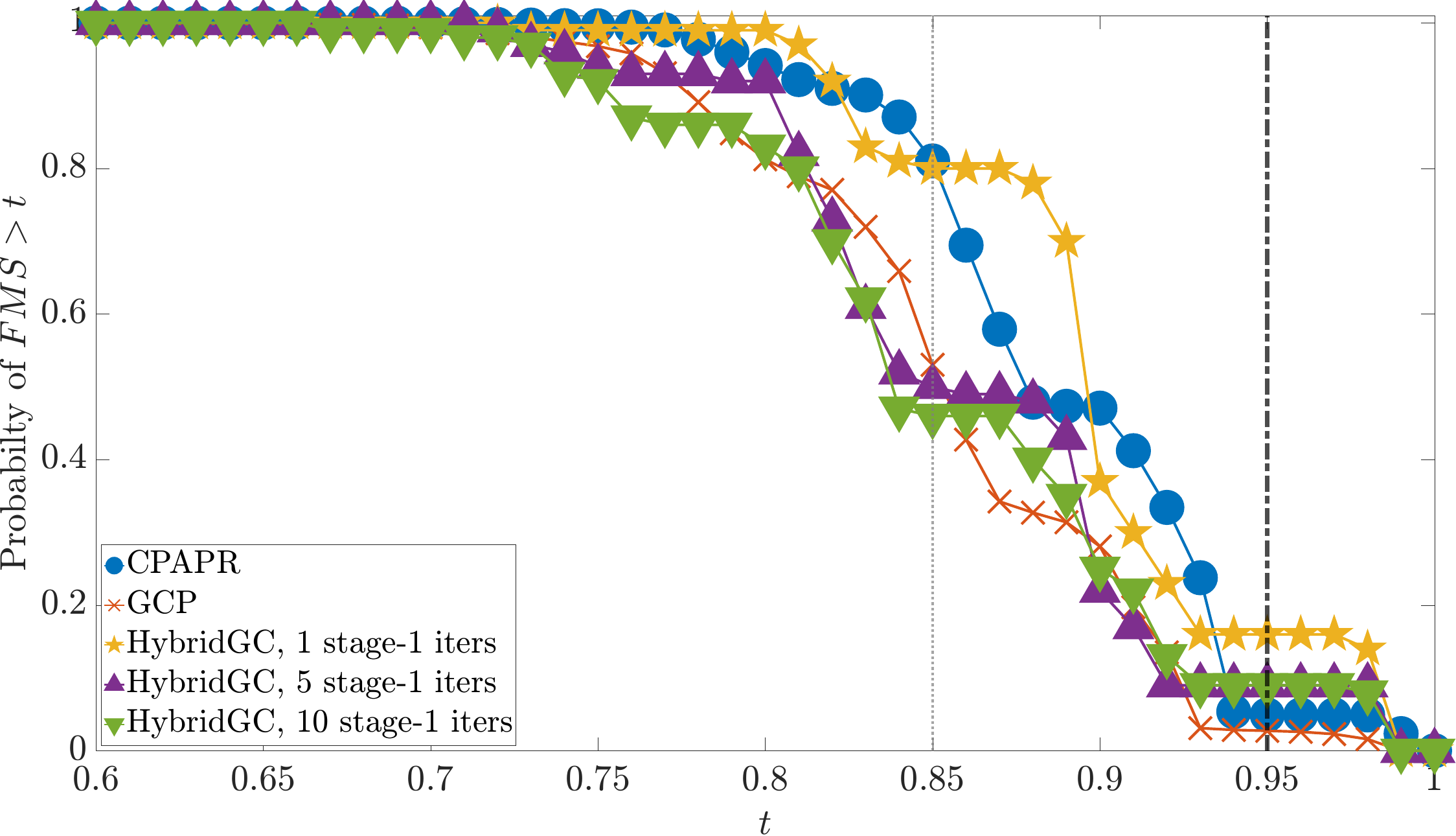}
	\caption{Factor match scores between CP models computed with
		\textsc{HybridGC}, CPAPR-MU, and GCP-Adam and the approximate global
		optimizer, $\tns{\widehat{M}}_{\mathcal{S}}^*$. The dash-dot gray
		vertical	lines and dotted black vertical lines denote the levels of
		``similar'' and ``equal'' described
		in~\cite{Lorenzo-Seva06TuckerCongruenceCoefficient}.}
	\label{fig:hybrid-gc:experiments:fms:prob_mle_est:1}
\end{figure}

Next, we evaluate \textsc{HybridGC} as a method for computing an approximate
low-rank basis to the global optimizer. We calculated the fraction of trials
with $\FMS$ greater than $t$~\cref{eq:hybrid-gc:error:fms:prob_mle_est} for GCP
and CPAPR, i.e., $\Psi(\tns{\widehat{M}}_{\mathcal{S}}^*,\mathcal{S}_G,t)$ and
$\Psi(\tns{\widehat{M}}_{\mathcal{S}}^*,\mathcal{S}_C,t)$, with $t \in [0,1]$.
We repeated this calculation for \textsc{HybridGC}, i.e.,
$\Psi(\tns{\widehat{M}}_{\mathcal{S}}^*,\mathcal{S}_H,t)$, and grouped the
results by the number of epochs taken by the first stage of
\textsc{HybridGC}. \Cref{fig:hybrid-gc:experiments:fms:prob_mle_est:1} presents
these results for GCP, CPAPR, and \textsc{HybridGC} (up to 10 epochs
of GCP). See \cref{fig:hybrid-gc:experiments:fms:prob_mle_est:2-10} in
\Cref{sec:hybrid-gc:experiments:supp} for supplementary results. Since all
curves showed the same behavior for $t < 0.6$, we report values for $t \in
	[0.6,1]$.

\textsc{HybridGC} tended to have a higher likelihood than GCP or CPAPR in
finding a low-rank basis equal to the empirical MLE when FMS $>0.95$, which is
considered high accuracy. This figure provides numerical evidence that
\textsc{HybridGC}---parameterized as a small amount of stochastic search ($\sim
	10$ GCP epochs) followed by deterministic search---was superior to GCP and
CPAPR by themselves in computing high accuracy models (solutions with FMS
greater than 0.95).

\section{Restarted CPAPR with \textsc{SVDrop}}
\label{sec:cpapr-svdrop}
\CHANGED{%
  Now, we present our second algorithm, Restarted CPAPR with \textsc{SVDrop}.%
}

\CHANGED{%
  We observed in~\cref{sec:hybrid-gc:experiments:nll}%
}
that there are situations where
\textsc{HybridGC} converges to the MLE but a standalone method like CPAPR does
not, and vice versa. The standalone method may converge to a minimizer far from
the MLE despite having started from the same point. The standalone method may
also converge to the MLE whereas \textsc{HybridGC} converges to some other
minimizer. Thus it is necessary to characterize the situations that end in
algorithm failure\footnote{ By ``success'', we mean convergence to the MLE. By
	``failure'' we mean convergence to any other KKT point \underline{or} failure to
	converge to a KKT point.} so that we may explain these seemingly conflicted
outcomes. Since it is easier to reason about a deterministic search path, we
focus on sources of failure in CPAPR and leave a similar study of GCP for future
work. CPAPR is also interesting because of its high success rate, meaning it is
sometimes feasible to examine all failed trials exhaustively. Furthermore only a
small number of parameters affect the search path of CPAPR in~\cite[Alg.
	3]{Chi12TensorsSparsityNonnegative}. Another motivator for the work presented
here is that CPAPR is used to refine the solution in \textsc{HybridGC}, thus
understanding convergence properties is important. Taking these factors into
account, CPAPR is a good candidate to analyze in order to understand why and
when local methods fail for Poisson CPD.

\CHANGED{%
  In this section we develop a heuristic called \textsc{SVDrop} to identify iterates of CPAPR
  that will converge to local minimizers based on the singular values of the tensor
  unfoldings. This heuristic is combined with a restarting procedure to
  form a novel variant of the CPAPR algorithm called Restarted CPAPR with
  \textsc{SVDrop} (\cref{alg:cpapr_svdrop}).%
}
We first motivate our analysis by observing a previously unreported problem and
reasoning as to its implications in \cref{sec:cpapr-svdrop:problem-context}.
In \cref{sec:cpapr-svdrop:rank-deficient}, we characterize this problem and
demonstrate empirically that: 1) the problem occurs frequently when computing
the Poisson CPD for a small synthetic dataset and 2) it can be
identified by tracking the $R$-th largest singular values of the mode
unfoldings at successive iterations.
The method is fully described in~\cref{sec:cpapr-svdrop:method}.
In
\cref{sec:cpapr-svdrop:experiments}, we present experimental evidence that
Restarted CPAPR with \textsc{SVDrop} improves the probability of convergence to
the MLE with an acceptable increase in computational cost.

\subsection{The drawback of extra (or too few) inner iterations}
\label{sec:cpapr-svdrop:problem-context}
Chi and Kolda's CPAPR paper~\cite{Chi12TensorsSparsityNonnegative} included a
section titled ``the benefit of extra inner iterations''. Their conclusion was
that although the maximum allowable number of inner iterations $l_{max}$ ``does
not significantly impact accuracy... increasing $l_{max}$ can decrease the
overall work and runtime''. They drew their conclusion from the mean and median
factor match score (FMS) between the model and the ``true solution''. However,
this definition of accuracy is incomplete since it ignores error estimators on
the loss function. Instead, our definition of accuracy includes both the
objective function value and expected convergence behavior over many trials.
Ultimately, we reach the opposite conclusion: the maximum allowable number of
inner iterations can significantly impact algorithm accuracy. Subsequently,
overall work and runtime are also affected. For instance, we observed situations
where, from one fixed starting point, CPAPR would  converge to different minima
for increasing values of $l_{max}$ in an alternating fashion: to the MLE for
some value of $l_{max}$, then to a different minimizer for a larger value, and
again to the MLE for an even larger value of $l_{max}$. In some cases, this
alternating pattern repeated multiple times. This behavior was not rare. In one
trial from 3,677 starting points, we observed some type of alternating
convergence pattern in 3,180 instances (86.4\%). In general, characterizing the
sensitivity of CPAPR to $l_{max}$ is complicated.

\begin{figure}[!ht]
	\centering
	\includegraphics[width=\textwidth]{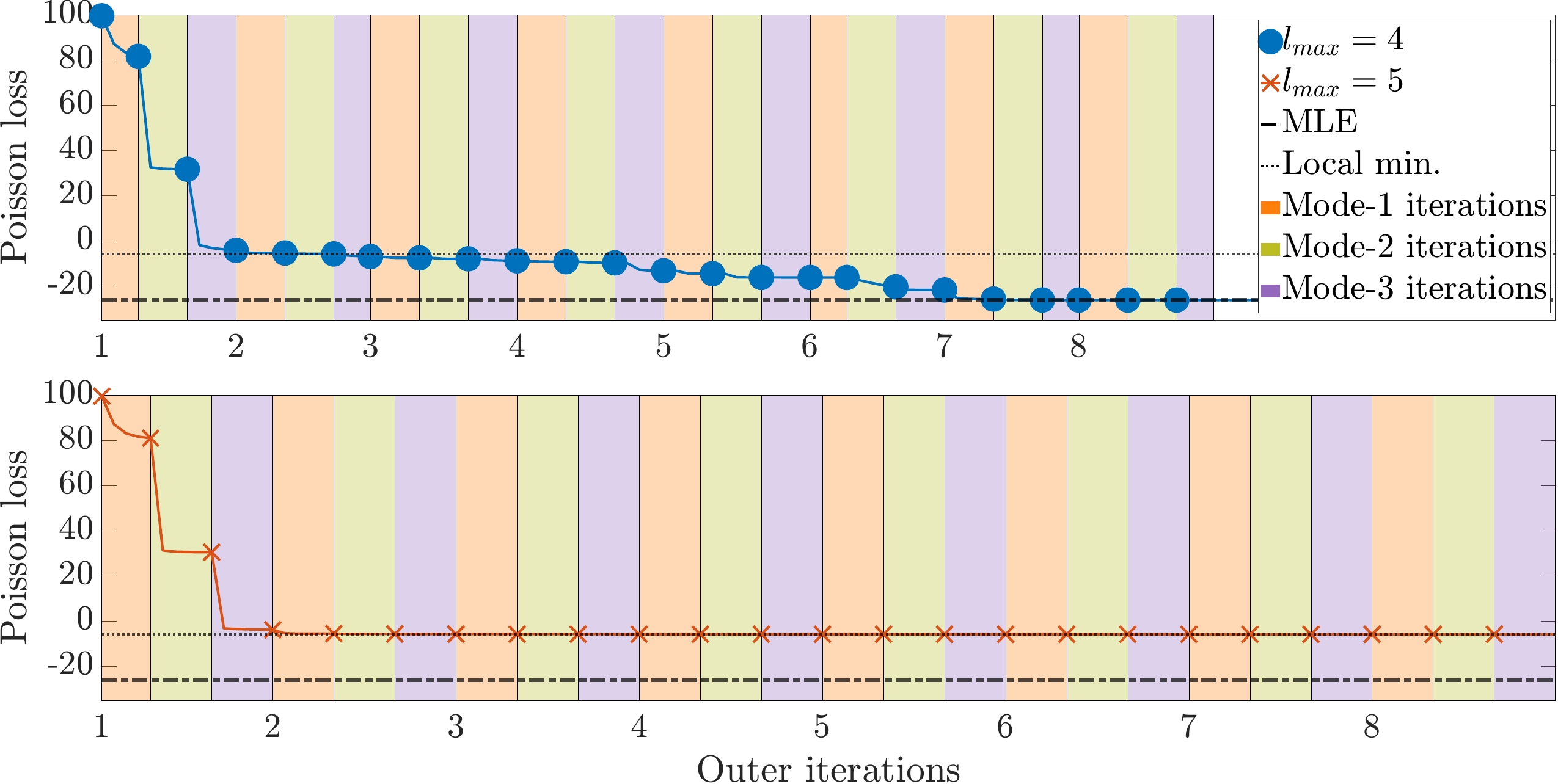}
	\caption{Traces of objective function values for the \emph{exemplar trial}:
    two decompositions computed by CPAPR starting from the same
    initial guess but with different numbers of maximum allowable
    inner iterations per mode.
    The $x$-axis
    is given in terms of the number of outer iterations, so
    optimizations by mode are differentiated with vertical blocks of
    color. Only the first 8 outer iterations are shown: CPAPR with
    $l_{max}=4$ (top) converges to the MLE; CPAPR with $l_{max}=5$
    (bottom) has settled in the basin of attraction of a different
    minimizer. The first Poisson loss value in each mode is emphasized
    with a marker.}
    \label{fig:cpapr-svdrop:inneriters-loop:example}
\end{figure}

\Cref{fig:cpapr-svdrop:inneriters-loop:example} demonstrates one example
empirically where the effect of extra inner iterations counters Chi and Kolda's
claim.
\CHANGED{%
	Both plots show the trace of the objective function value for the first 8 outer iterations of
	CPAPR.
	The upper plot shows the trace of the objective function value over the
	total iteration history when $l_{max}=4$.
	CPAPR converges to the empirical MLE (black
	dash-dot line) in 8 outer iterations.
	The lower plot is taken from the same initial
	point except the maximum number of inner iterations is one greater, i.e., $l_{max}=5$.
	By the 8th outer iteration, CPAPR has settled in the basin of
  attraction of a KKT point far from the MLE (black dotted line).%
}
We will return to this case
throughout the rest of this section, so we will refer to it as the
\emph{exemplar trial}.

\begin{figure}[!ht]
	\centering
	\includegraphics[width=\textwidth]{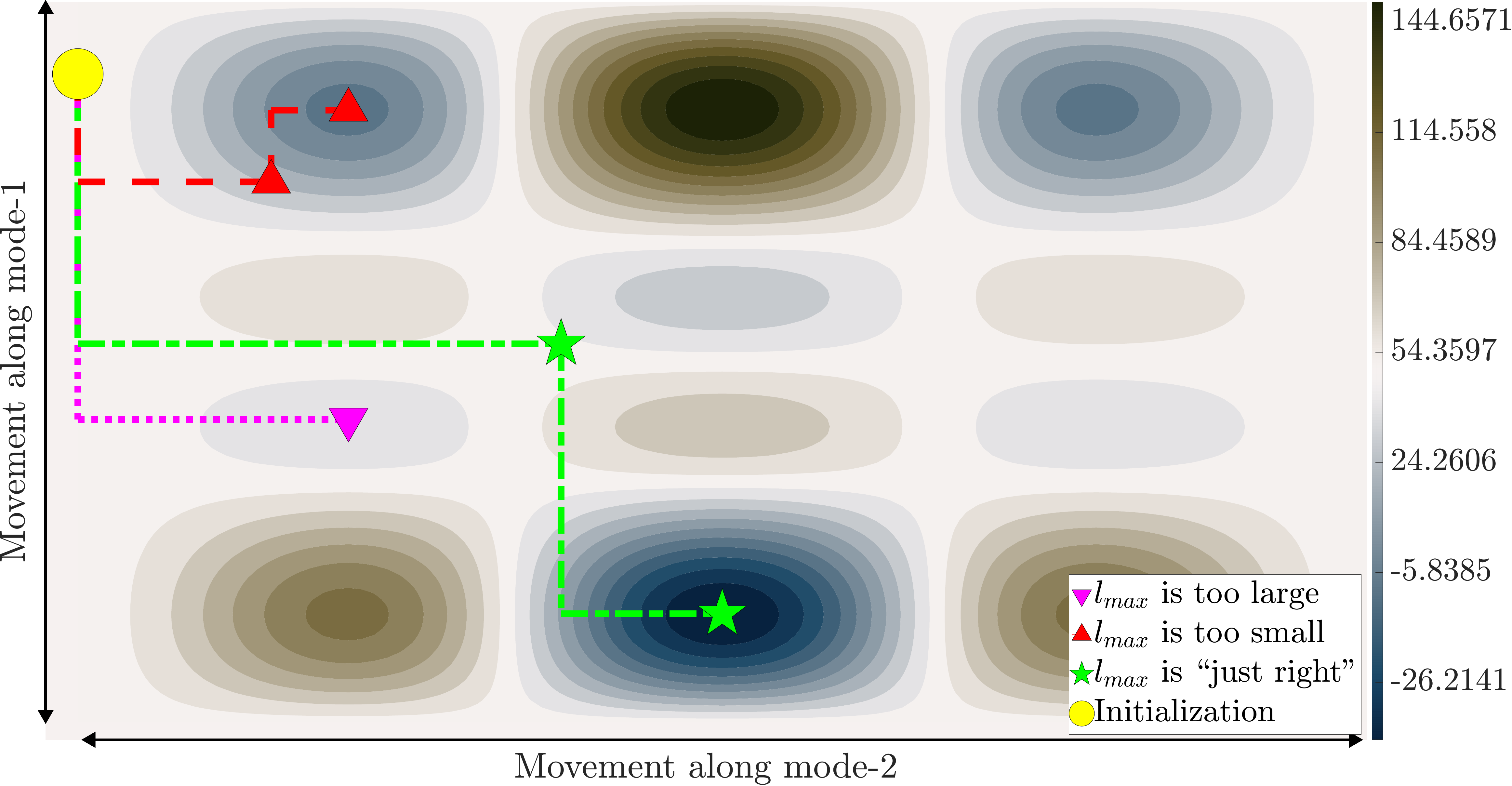}
	\caption{The contour plot illustrates how convergence depends on the number
		of inner iterations in the search direction for a 2D problem. Blue
		represents minima and brown represents maxima; darker shades are more
		extreme values than lighter shades. From the same starting
		initialization, CPAPR is run with three different values for inner
		iterations $l_{max}$: 1) the ``Goldilocks'' amount that leads to the
		MLE; 2) too few or 3) too many inner iterations, which both lead to
		different minimizers.}
	\label{fig:cpapr-svdrop:inneriters-loop:idea}
\end{figure}

Incremental changes to this parameter can result in drastically different
outcomes. In experimentation, we frequently observe that, for the first several
outer iterations, CPAPR tends to max out the number of inner iterations in each
mode without converging. This led to the conclusion that early iterations are
especially critical and sensitive to $l_{max}$.
\Cref{fig:cpapr-svdrop:inneriters-loop:idea} presents a conceptual model
explaining the situation. The contour plot reflects the minima (in blue) and
maxima (in brown) of a $d=2$ problem. Darker shades reflect more extreme values.
The $x$- and $y$-axes represents search paths in the directions of the second
and first modes, respectively. The green, magenta, and red lines represent the
search paths of CPAPR from the same initial starting point (yellow circle) but
with different values for $l_{max}$. The magenta and red paths allow too many or
too few inner iterations and converge to local minimizers. The green path allows
the ``Goldilocks'' amount---this choice leads to the MLE.\footnote{The fairy
	tale of Goldilocks is about a girl named Goldilocks who enters a house belonging
	to three bears and tries out their belongings to find one that suits her best.
	Each item belonging to one of the three bears that she samples is either ``too
	many'', ``too few'', or ``just right''. Allusions to the Goldilocks fairy tale
	are also used by astronomers to describe the
	\href{https://exoplanets.nasa.gov/resources/323/goldilocks-zone/}{zone of
		habitable exoplanets around a star}.} We will show that our novel analysis can
differentiate between the green path and the magenta and red paths at runtime.

\subsection{Spectral properties identify rank-deficient solutions}
\label{sec:cpapr-svdrop:rank-deficient}
\begin{figure}[!hp]
	\centering
	\begin{subfigure}[b]{1.00\textwidth}
		\centering
		\includegraphics[width=\textwidth]{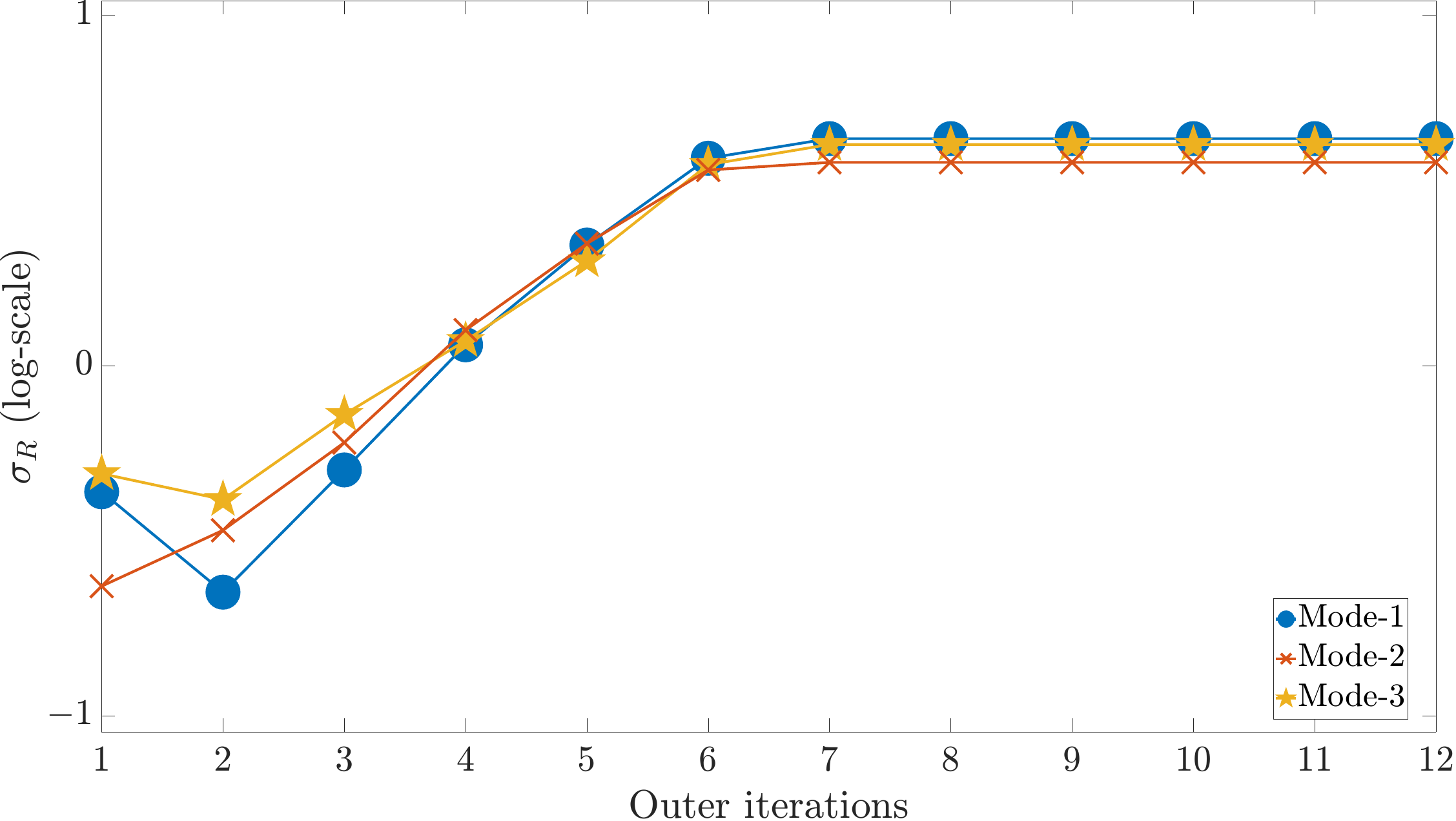}
		\caption{$l_{max} = 4$.}
		\label{fig:hybrid-gc:experiments:sigmas:traces:mle}
	\end{subfigure}\\
	\begin{subfigure}[b]{1.00\textwidth}
		\centering
		\includegraphics[width=\textwidth]{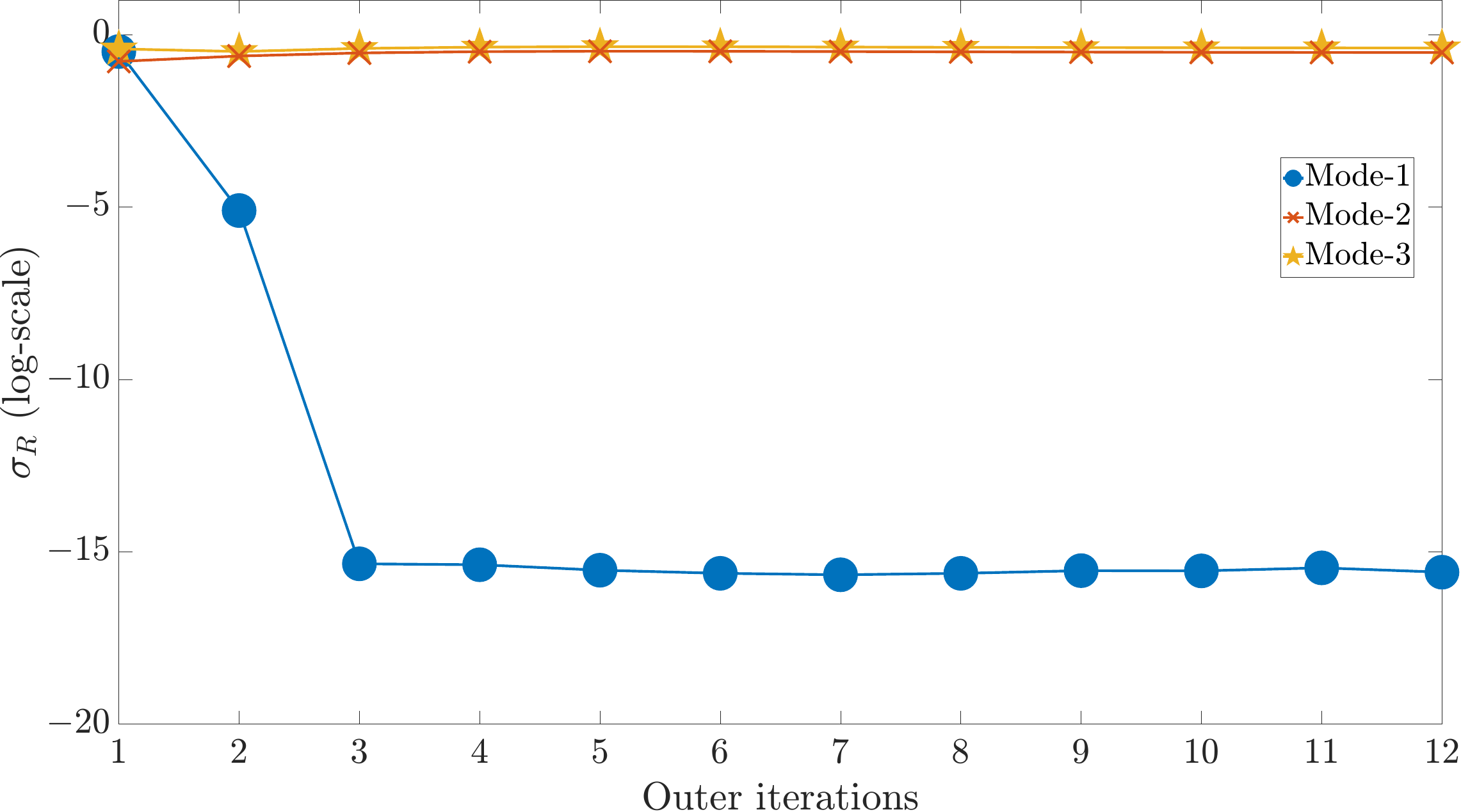}
		\caption{$l_{max} = 5$.}
		\label{fig:hybrid-gc:experiments:sigmas:traces:loc}
	\end{subfigure}
	\caption{\CHANGED{Traces of the $R=3$-rd largest singular value of each
		mode unfolding after every update for the
		\emph{exemplar trial} on \texttt{LowRankSmall}.
		Analyzing the $R$-th largest singular value of each
		mode yields an explicit signal indicating convergence
		to a rank-deficient solution.
		The $x$-axis in both
		plots has been truncated to the total number of outer
		iterations needed to converge to the MLE; the traces
		shown in the lower plot continue until convergence to
		the local minimizer without change.}
	}
	\label{fig:hybrid-gc:experiments:sigmas:traces}
\end{figure}

By analyzing the singular values of each mode unfolding, we can see critical
changes to the model tensor that are otherwise hidden. In particular, we
consider the $R$-th largest singular value when the requested decomposition rank
is $R$. \Cref{fig:hybrid-gc:experiments:sigmas:traces} shows the values of the
third largest singular value (since $R=3$) of each of the mode unfoldings over
the iteration
history:~\cref{fig:hybrid-gc:experiments:sigmas:traces:mle} plots when
$l_{max} =4$, the case when CPAPR converges to the
MLE;~\cref{fig:hybrid-gc:experiments:sigmas:traces:loc} plots when
$l_{max} = 5$, the case when CPAPR converges to a different local
minimizer. Both trials iterated from the same initial guess. The critical
observation is that the $R$-th singular value in mode-1 when $l_{max}=5$ inner
iterations are taken per outer iteration (blue line with circle
markers in~\cref{fig:hybrid-gc:experiments:sigmas:traces:loc}) is driven below
machine precision. \CHANGED{We define a \emph{rank-deficient solution}
as an approximate solution (i.e., iterate of CPAPR) where one or more
of the $R$-th largest singular values of the tensor unfoldings across all modes is zero}. A rank-deficient
solution is a KKT point, but it is not the MLE. \emph{Ceteris paribus}, the
search path that leads to it is determined by $l_{max}$.

\begin{figure}[!hp]
	\centering
	\begin{subfigure}[b]{1.00\textwidth}
		\includegraphics[width=\textwidth]{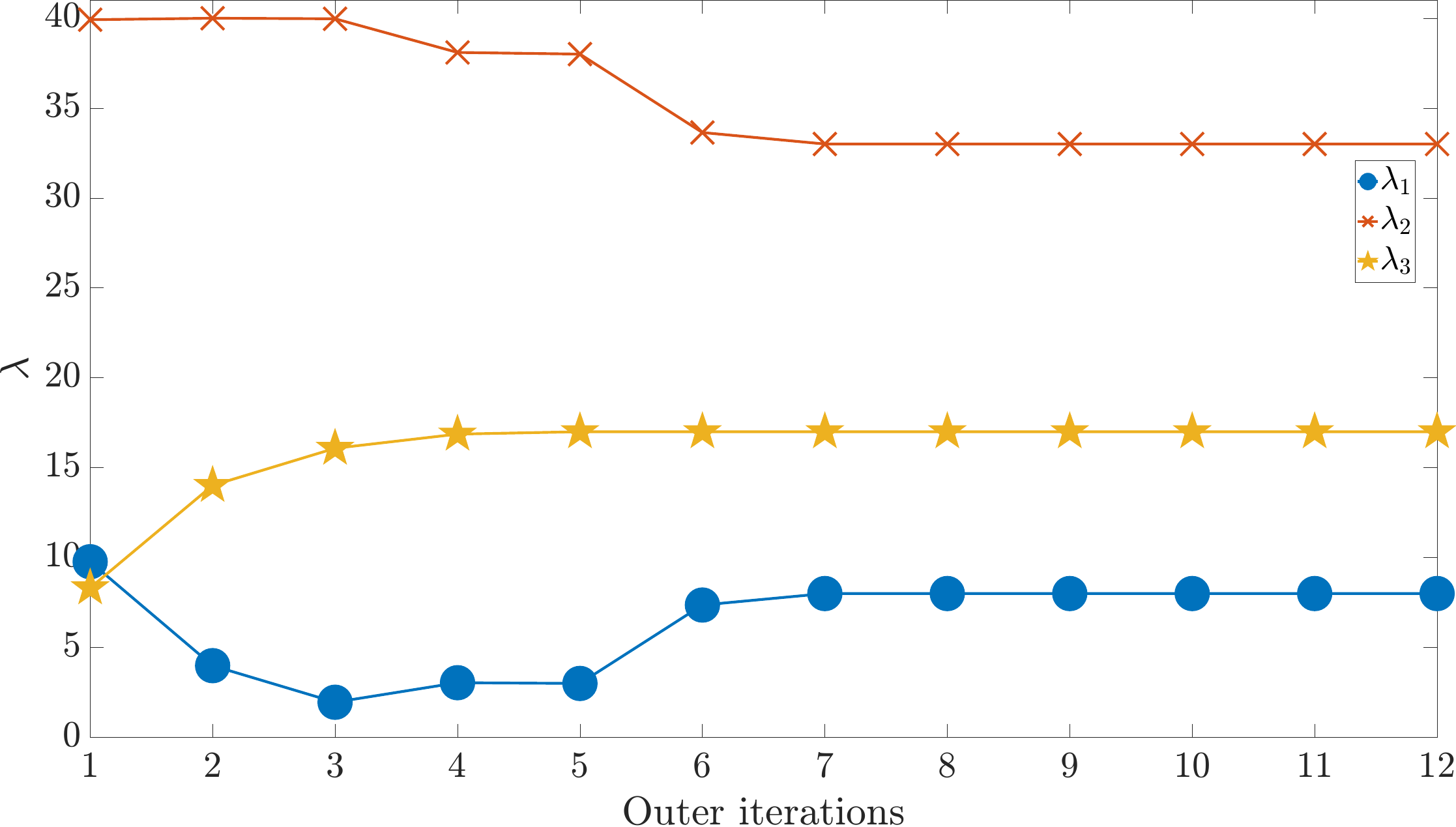}
		\caption{$l_{max}=4$.}
		\label{fig:hybrid-gc:experiments:lambdas:mle}
	\end{subfigure}\\
	\begin{subfigure}[b]{1.00\textwidth}
		\includegraphics[width=\textwidth]{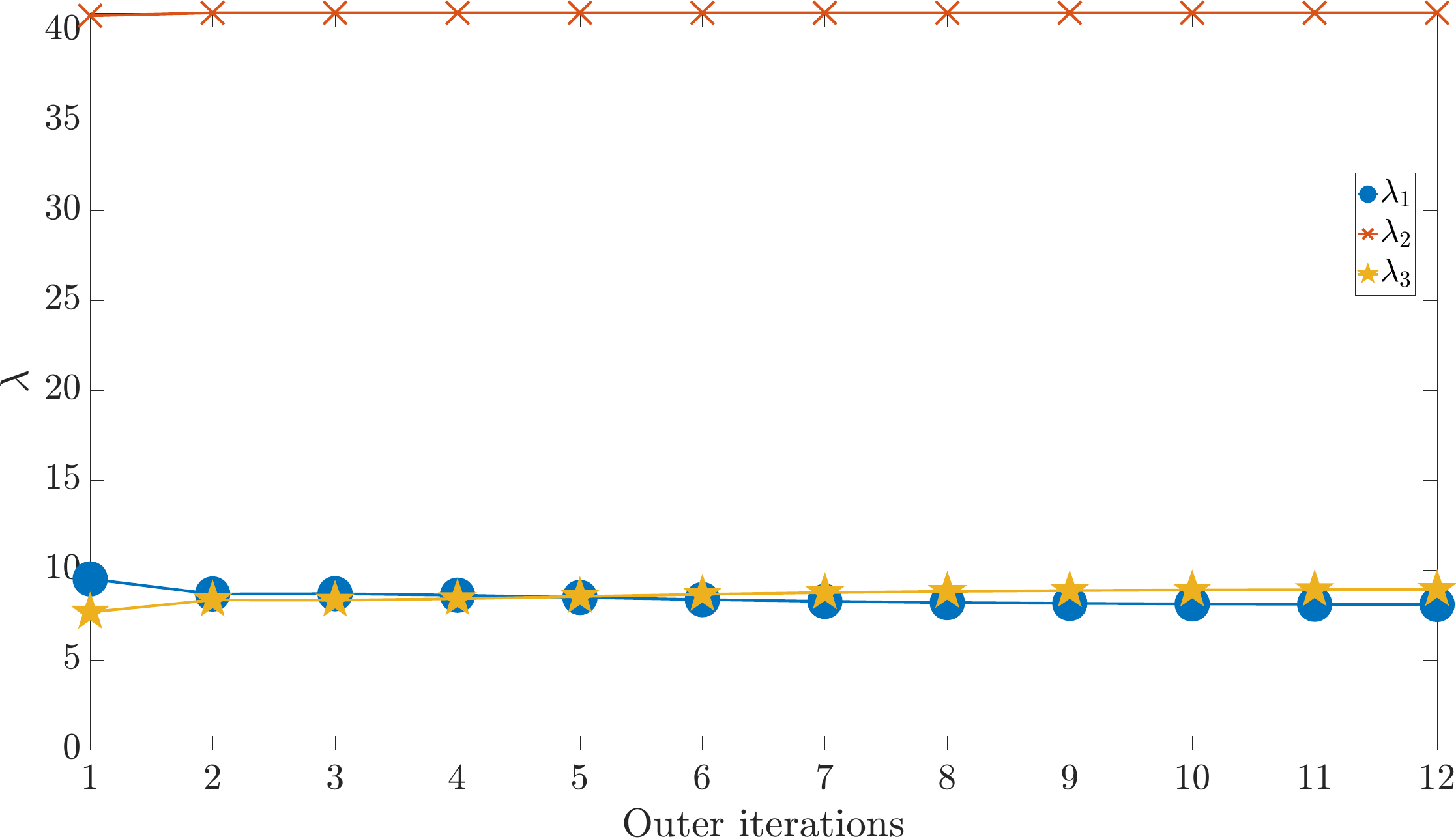}
		\caption{$l_{max}=5$.}
		\label{fig:hybrid-gc:experiments:lambdas:loc}
	\end{subfigure}
	\caption{\CHANGED{Traces of the $R$ $\lambda$ weights maintained by
		CPAPR in each iteration for the \emph{exemplar trial}
		on \texttt{LowRankSmall}. It is unclear whether the
		$\lambda$-weights provide a useful heuristic for determining
		convergence to the MLE versus a rank-deficient solution.
		The $x$-axis in both
		plots has been truncated to the total number of outer
		iterations needed to converge to the MLE; the traces
		shown in the lower plot continue until convergence to
		the local minimizer without change.}
	}
	\label{fig:hybrid-gc:experiments:lambdas:traces}
\end{figure}

At present this analysis of success versus failure is only possible by examining
the singular values of the mode unfoldings.
\CHANGED{%
  The $\lambda$ values of the CP model, which are freely
  available, do not provide the same information as the singular
  values of the mode unfoldings, which must be computed at additional cost.
  For example,~\cref{fig:hybrid-gc:experiments:lambdas:traces}
}
displays the $\lambda$ values of the CP model at each iteration for
both values of $l_{max}$ in the \emph{exemplar trial}. The $\lambda_1$
and $\lambda_3$ values (blue and yellow, respectively) became nearly
identical when the number of inner iterations taken per outer
iteration was $l_{max}=5$. This occurred at nearly the same time the
mode-1 $R$-th singular value started to drop drastically. It remains
unclear whether this behavior may indicate convergence to
a rank-deficient solution or if it is just a coincidence. We did not
observe a similar pattern in the $\lambda$ weights in a sample of
other trials with similar behaviors in the singular values. The
behavior also is not explained by standard techniques in optimization,
e.g., unsatisfied convergence criteria, or machine learning, e.g.,
large gradient norms. The $R$-th largest singular value of the mode-1
unfolding is a pronounced indicator of convergence to a rank-deficient
solution.


We observed this behavior in most cases. Recall from
\cref{tab:hybrid-gc:experiments:nll:prob_mle:data1} that CPAPR converged to the
MLE in 106,215 of 110,266 trials ($\widehat{P}(\epsilon = 10^{-4}) = 0.963$).
For a random sample of 10,000 of these trials,\footnote{Computing this statistic
	for all random starts was prohibitively expensive.} the Poisson CP models were
not rank-deficient: the $R$-th largest singular value when CPAPR terminated was
typically far from machine precision ($4.449$ on average). By contrast, in 4,020
of the 4,051 trials (99.235\%) where CPAPR converged to a different KKT point,
we found that the $R$-th largest singular value in mode-1 when CPAPR terminated
was on the order of double precision machine epsilon (i.e., $\approx 2.2204
	\times 10^{-16}$). Thus we describe these models as rank-deficient: when the
$R$-th largest singular value in one or more modes is numerically close to 0
(i.e., near or below machine precision) so that the column rank of these mode
unfoldings is less than $R$. We can reasonably conclude that there is a strong
connection between rank-deficient solutions and KKT points that are not the
MLE.

\subsection{Restarted CPAPR with the \textsc{SVDrop} heuristic}
\label{sec:cpapr-svdrop:method}
Upon closer examination of the traces from the rank-deficient search path
in~\cref{fig:hybrid-gc:experiments:sigmas:traces}, we notice that the $R$-th
largest singular value in mode-1 drops dramatically between the 29th and 30th
iteration. (In this case, $l_{max}=	5$, so the 30th total iteration is the first
inner iteration on mode-1 of the 3rd outer iteration.) To be precise, the
\emph{gap ratio} of the $R$-th largest singular value from the mode-1 unfolding
between iterations 29 and 30 is $\sigma_{(1)}[R]^{(29)}/\sigma_{(1)}[R]^{(30)}
	\approx 3.6\times 10^{10}$. Considering all of the failed trials, the maximum
gap ratio was $1.55\times 10^{12}$ on average, the median of the maximum gap
ratios was $2.95 \times 10^6$, and the median iterate where it was observed was
the 30th iteration. Analogous to the indication of numerical instability by
large condition number, we interpret a large gap ratio between successive
iterates as indicative of a search path that will converge to a rank-deficient
solution. Therefore, large gap ratio may serve as a reliable heuristic to
determine whether the current search path should be accepted or rejected. This
observation informs our solution in \Cref{alg:cpapr_svdrop}. Before we
present the method in full, two core concepts remain to be discussed: 1) the
\textsc{SVDrop} heuristic and 2) Restarted CPAPR.

\CHANGED{%
  In our experiments, we use the MATLAB dense \texttt{svd} solver to compute
	gap ratios. Regardless of the SVD solver chosen, calculating the gap ratio between
	successive iterates will add non-trivial costs due to the SVD
	computation.\footnote{When $\mat{A}$ is dense, the fastest \emph{direct} or
		\emph{transformation method} for serial or shared-memory parallel execution to
		compute all singular triplets of the matrix incurs a cost of $4mn^2
			- \frac{4}{3}n^3$ floating point operations (FLOPS) or
		better~\cite{Bai00TemplatesSolutionAlgebraic,
			Golub65CalculatingSingularValues, Golub96MatrixComputations,
			Parlett80SymmetricEigenvalueProblem} (assuming $m\geq n$). It is accessible
		through LAPACK's \texttt{xGESVD} driver~\cite{Anderson99LAPACK}.} Furthermore,
	it is not even clear whether the gap ratio needs to be computed after every
	update. To mitigate incurring these excessive computational costs, we propose the
	\textsc{SVDrop} heuristic in a new subproblem for an extended CPAPR algorithm
	in~\cref{sec:cpapr-svdrop:procedure}.
	Additionally, in~\cref{sec:cpapr-svdrop:discussion:svd} we discuss how iterative SVD
	solvers (e.g., Krylov methods or randomized SVD) may be used to further reduce
	costs.%
}

\subsubsection{Procedure}
\label{sec:cpapr-svdrop:procedure}
The inputs to the subproblem algorithm are:
\begin{enumerate}
	\item the maximum number of inner iterations, $l_{max} \in \N_{>0}$;
	\item the number of \textsc{SVDrop} inner iterations between successive
	      models for computation of their spectral properties, $\tau \in \Z_+$;
	      and,
	\item a maximum threshold for the gap ratio indicating an acceptable search
	      path, $\gamma \in \R_{+}$.
\end{enumerate}

\CHANGED{%
	While the model is not converged in~\Cref{alg:cpapr_svdrop}, we compute
	a rank-$R$ decomposition with CPAPR and improve the model fit with
	multiplicative updates and calls to~\Cref{alg:cpapr_svdrop:subproblem}.
	The CPAPR subproblem loop begins at
	Line~\ref{alg:cpapr_svdrop:subproblem:subproblem_loop}
	in~\Cref{alg:cpapr_svdrop:subproblem}. Every
	$\tau$ inner iterations, we compute the gap ratio between the current model and
	the checkpointed model from the singular values of the mode unfoldings at
	Lines~\ref{alg:cpapr_svdrop:subproblem:svd_iter_l}--\ref{alg:cpapr_svdrop:subproblem:compute_gamma}.
	If the gap ratio computed is smaller than $\gamma$, then we accept the search
	path, checkpoint the model, and continue iterating. When the
	gap ratio is greater than $\gamma$, the \textsc{SVDrop} heuristic indicates that the search path
	will converge to a
	rank-deficient solution and we reject it, since to otherwise continue likely
	will waste a great deal of computation, as seen in
	\cref{fig:hybrid-gc:experiments:nll:data1:traces} and
	\cref{fig:hybrid-gc:experiments:sigmas:traces}. A simple option is to
	restart: discard the work done up to now, randomly choose a new starting point
	in the feasible domain of the optimization problem, and recompute. Restarting,
	taken together with the \textsc{SVDrop} heuristic, is the idea behind our new
	method, Restarted CPAPR with \textsc{SVDrop}, presented in
	\Cref{alg:cpapr_svdrop} and \Cref{alg:cpapr_svdrop:subproblem}.
	The exposition follows the template of the ideal version of CPAPR (see~\cite[Algs.
		1--2]{Chi12TensorsSparsityNonnegative}).
	The symbol $\odot$ denotes the \emph{matricized tensor times Khatri-Rao product
		(MTTKRP)}~\cite{Smilde04MultiwayAnalysisApplications}, an important computational
	kernel in many tensor decomposition algorithms. The symbol $\oslash$ denotes elementwise
	division. The symbol $\ast$ denotes elementwise matrix multiplication, i.e., the
	Hadamard product. We denote $\mat{e}$ as a vector of all ones with the
	appropriate dimension corresponding to mode $n$
	in~\Cref{alg:cpapr_svdrop} and~\Cref{alg:cpapr_svdrop:subproblem},
	$\mat{\lambda}\in\R^d$ as the vector storing the weights of the Kruskal tensor
	$\tns{M}$ in~\Cref{eq:cp}, and $\diag(\mat{\lambda}) \in\R^{d \times d}$
  is a square matrix with the elements of $\mat{\lambda}$ along the
  diagonal and zeros everywhere
	else.%
}

\subsubsection{Additional considerations} First, we reset the counter of
\textsc{SVDrop} inner iterations before the first inner iteration of a mode,
since we observed that  the singular values of the unfoldings of the modes that
are held fixed change very little between iterations.\footnote{It is possible
	that the singular values of the mode unfoldings held fixed may change. However,
	this was beyond the scope of this work.} A consequence is that the number of
\textsc{SVDrop} inner iterations should always be less than or equal to the
maximum number of inner iterations, i.e., $\tau \leq l_{max}$. Second, the gap
ratio tolerance $\gamma$ should be sufficiently large. We used $\gamma = 10^6$
in our numerical experiments but it is unknown if this value generalizes to
other problems. Third, the update step Line~\ref{alg:cpapr_svdrop:subproblem:updateMwithB}
in \Cref{alg:cpapr_svdrop:subproblem} is for exposition; it can be
performed implicitly on $\tns{M}$ efficiently in software.

\begin{algorithm}[!ht]
	\caption{%
		Restarted CPAPR algorithm with \textsc{SVDropSubproblem}.%
	} 
	\label{alg:cpapr_svdrop}
	\begin{algorithmic}[1]
		\Statex \noindent Let $\tns{X}$ be a tensor of size $I_1 \times \dots
			\times I_d$.
		\Statex User input:
		\begin{itemize}
			\item $R$: Number of components in low-rank approximation
			\item $l_{max}$: Maximum number of inner iterations per outer
			      iteration
			\item $\tau$: Number of \textsc{SVDrop} inner iterations to perform
			\item $\gamma$: Tolerance for identifying rank-deficiencies (e.g.,
			      $10^6$)
		\end{itemize}
		\State $\texttt{restart} \gets \text{true}$
		\Repeat \, \LeftComment{\textsc{SVDrop}}
		\If{$(\texttt{restart})$}
		\State $\tns{M} \gets \Call{GenerateRandomGuess}{[I_1,\ldots,I_d],R}$
		\label{alg:cpapr_svdrop:restart}
		\EndIf
		\Repeat \, \LeftComment{\textsc{OuterIteration}}
		\label{alg:cpapr_svdrop:outer_loop}
		\For{$n = 1,\ldots,d$}
		\label{alg:cpapr_svdrop:mode_loop}
		\State $\mat{\Pi} \gets \left( \mat{A}^{(d)} \odot \dots \odot \mat{A}^{(n+1)} \odot \mat{A}^{(n-1)} \odot \dots \odot \mat{A}^{(1)} \right)^T$
		\label{alg:cpapr_svdrop:mttkrp}
		\State $\left[ \mat{B}, \texttt{restart} \right] \gets \Call{SVDropSubproblem}{\tns{M},\mat{\Pi},R,n,l_{max},\tau,\gamma}$
		\label{alg:cpapr_svdrop:call_subproblem}
		\If{$(\texttt{restart})$}
		\State \textbf{break}
		\label{alg:cpapr_svdrop:call_restart}
		\Comment{Break out of \{\textsc{OuterIteration}\}.}
		\EndIf
    \State
    \CHANGED{%
        $\mat{\lambda} \gets \mat{e}^T\mat{B}, \,\mat{A}^{(n)} \gets \mat{B} \diag(\mat{\lambda})^{-1}$
    }
		\label{alg:cpapr_svdrop:updateMwithB}
		\EndFor
		\Until convergence
		\LeftComment{\textsc{OuterIteration}}
		\Until convergence
		\LeftComment{\textsc{SVDrop}}
	\end{algorithmic}
\end{algorithm}

\begin{algorithm}[!ht]
	\caption{%
		Subproblem solver for Algorithm~\ref{alg:cpapr_svdrop} with
		\textsc{SVDrop} heuristic.%
	} 
	\label{alg:cpapr_svdrop:subproblem}
	\begin{algorithmic}[1]
		\Function{SVDropSubproblem}{\CHANGED{%
				Kruskal tensor $\tns{M} = \llbracket
				\mat{\lambda};\,\mat{A}^{(1)}, \ldots, \mat{A}^{(d)}
				\rrbracket$,
			a sequence of Khatri-Rao products $\mat{\Pi}$, desired rank $r$, mode $n$, maximum number
			of inner iterations $l_{max}$, number of \textsc{SVDrop} inner iterations $\tau$,
			maximum threshold for gap ratio $\gamma$%
      }
    }
		\State $t \gets 0$ \State $\sigma_r \gets \text{svd}(\mat{M}_{(n)})[r]$
		\label{alg:cpapr_svdrop:subproblem:svd_iter_0}
		\Comment{$r$-th largest singular value of $\mat{M}_{(n)}$.}
		\State \CHANGED{%
      $\mat{B} \gets \mat{A}^{(n)} \diag(\mat{\lambda})$
    }
		\label{alg:cpapr_svdrop:subproblem:distribute}
		\For{$l = 1,\ldots,l_{max}$}
		\label{alg:cpapr_svdrop:subproblem:subproblem_loop}
		\State $t \gets t + 1$
		\State $\mat{\Phi} \gets \left( \mat{X}_{(n)} \oslash \left(\mat{B}\mat{\Pi}\right)\right)\mat{\Pi}^T$
		\label{alg:cpapr_svdrop:subproblem:compute_phi}
		\State $\mat{B} \gets \mat{B} \ast \mat{\Phi}$
		\label{alg:cpapr_svdrop:subproblem:update_B_with_Phi}
		\If{$(t = \tau)$}
		\State \CHANGED{%
      $\mat{\lambda} \gets \mat{e}^T\mat{B}, \,\mat{A}^{(n)} \gets
      \mat{B} \diag(\mat{\lambda})^{-1}$%
    }
		\label{alg:cpapr_svdrop:subproblem:updateMwithB}
		\Comment{Update $\mat{A}^{(n)}$, which updates $\tns{M}$.}
		\State $\sigma_r^{(l)} \gets \text{svd}(\mat{M}_{(n)})[r]$
		\label{alg:cpapr_svdrop:subproblem:svd_iter_l}
		\If{$(\sigma_r / \sigma_r^{(l)} > \gamma)$}
		\label{alg:cpapr_svdrop:subproblem:compute_gamma}
		\State \Return $[\mat{B}, \texttt{true}]$
		\label{alg:cpapr_svdrop:subproblem:force_restart}
		\Comment{%
      \CHANGED{%
        $\rank(\mat{M}_{(k)})<r$;%
      }
		forces a restart in~\Cref{alg:cpapr_svdrop}.}
		\EndIf
		\State $\sigma_r \gets \sigma_r^{(l)}$
		\State $t \gets 0$
		\EndIf
		\EndFor
		\State \Return $[\mat{B},\texttt{false}]$
		\Comment{%
      \CHANGED{%
        $\rank(\mat{M}_{(k)})=r$;%
      }
		does not force a restart in~\Cref{alg:cpapr_svdrop}.}
		\EndFunction
	\end{algorithmic}
\end{algorithm}

\subsection{Numerical experiments}
\label{sec:cpapr-svdrop:experiments}

To evaluate Restarted CPAPR with \textsc{SVDrop}, we selected 14,051 random
initializations for \texttt{LowRankSmall} from the experiments
in~\cref{sec:hybrid-gc:experiments}: the random sample of 10,000 starts where
CPAPR converged to the MLE that we mentioned previously plus the 4,051 starts
where CPAPR converged to a different KKT point. We computed rank $R=3$ CP
decompositions using CPAPR with Multiplicative Updates starting from each point.
We increased the number of \textsc{SVDrop} inner iterations as $\tau \in
\{0,\ldots,10\}$ but kept all other parameters identical to the experiments
in~\cref{sec:hybrid-gc:experiments}. A value $\tau = 0$ indicates that
\textsc{SVDrop} was not used and that CPAPR was not restarted at any iteration.

\subsubsection{Probability of convergence}

In the previous experiments using CPAPR without any restarts, the total number
of trials that did not converge to the MLE at the level of $\epsilon=10^{-4}$
was 4,051 of 110,266 ($1-\widehat{P}_{MLE} (10^{-4}) = 0.037$; see
\Cref{tab:hybrid-gc:experiments:nll:prob_mle:data1}).
\Cref{tab:cpapr-svdrop:results} reports the total number of trials that: 1)
converged to the MLE, 2) converged to some other KKT point, or 3) did not
converge to any KKT point using this set of 4,051 initial starts. Convergence
was calculated as in~\cref{eq:hybrid-gc:error:nll:prob_mle_est} at the level of
$\epsilon=10^{-4}$ for each value of $\tau$.\footnote{Our results hold to the
level of $\epsilon =10^{-8}$ but we report $\epsilon=10^{-4}$ to make comparison
with the previous experiments with \textsc{HybridGC} without any restarts.} Of
these, 3,905 converged to some other KKT point and 146 did not converge to any
KKT point when $\tau = 0$. Our method improved on this in all cases. It is
interesting to note that the probabilities of convergence to a different KKT
point and failure to converge to any KKT point were higher when $\tau \geq 6$.
We defer further discussion on this point to~\cref{sec:cpapr-svdrop:discussion}
since we will provide additional results that will help us reason about this
behavior and allow us to make suggestions with more context.

\begin{figure}[!ht]
	\begin{subfigure}[t]{1.00\textwidth}
		\includegraphics[width=\textwidth]{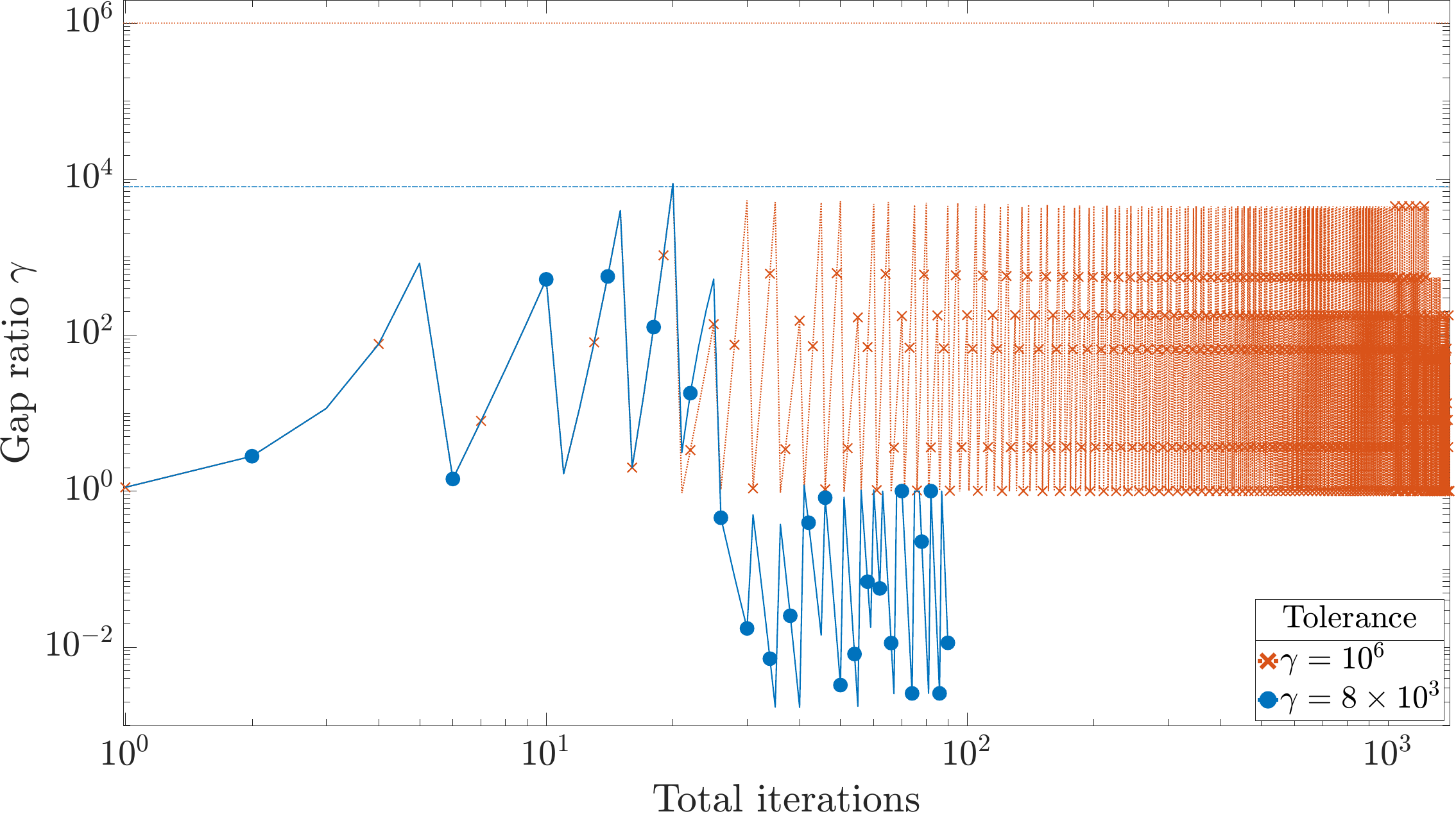}
	\end{subfigure}\\
	\begin{subfigure}[t]{1.00\textwidth}
		\includegraphics[width=\textwidth]{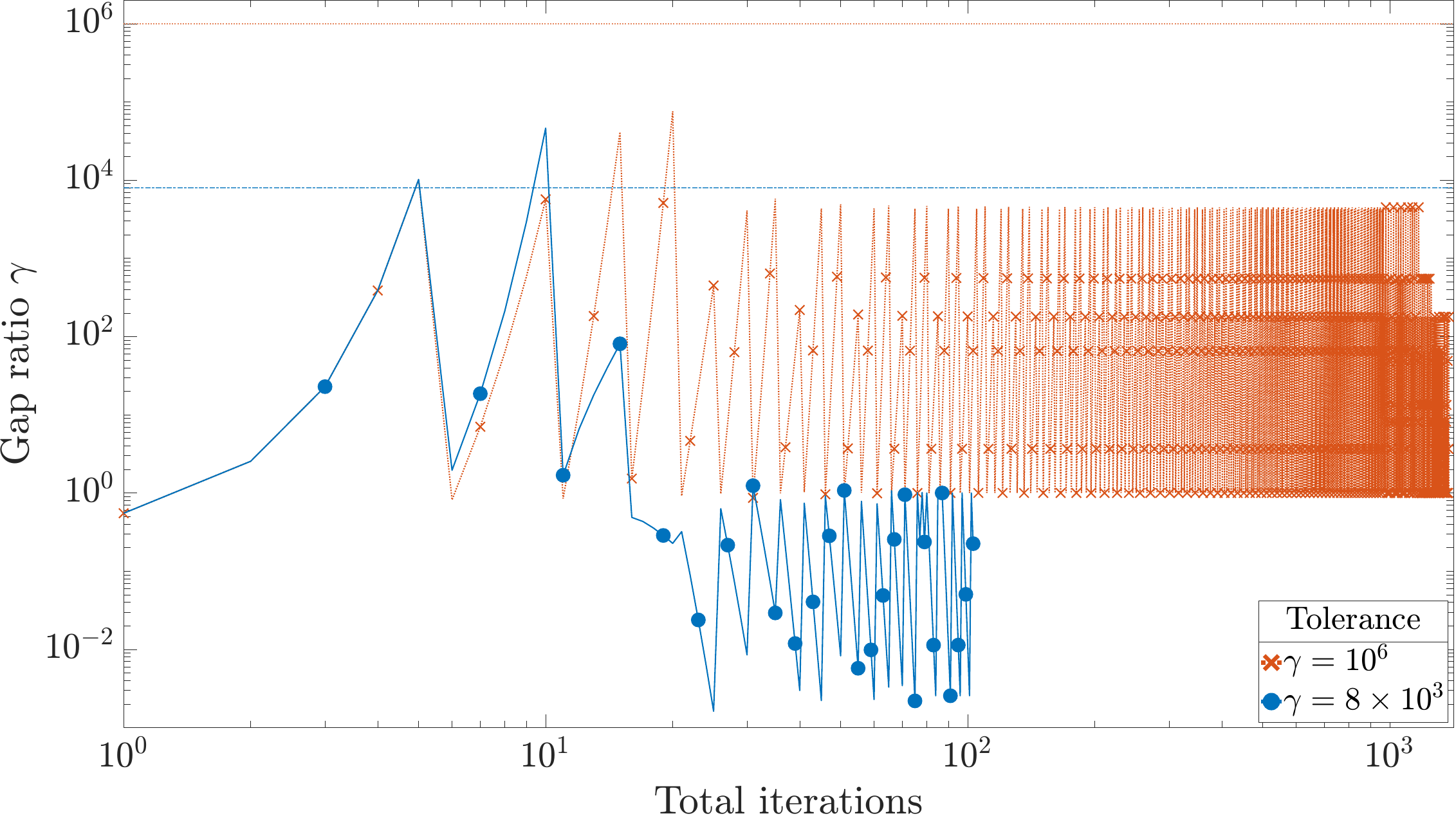}
	\end{subfigure}
	\caption{Two trials demonstrating the sensitivity of \textsc{SVDrop} to
		$\gamma$. Starting from points where CPAPR without restarting (i.e.,
		$\tau=0$) was known to converge to the MLE, \textsc{SVDrop} failed to
		converge to any KKT point when $\gamma$ was set too large (red
		$\times$ marker). When
		$\gamma$ is too large, \textsc{SVDrop} may fail to recognize rank
		deficiency and CPAPR may stagnate near a local minimizer that is not a
		KKT point. When set appropriately (blue circle marker), \textsc{SVDrop} converged to
		the MLE.}
	\label{fig:cpapr-svdrop:fail}
\end{figure}

In the best case, when the number of \textsc{SVDrop} inner iterations was
$\tau=2$, Restarted CPAPR with \textsc{SVDrop} recovered the MLE in all but two
trials ($\widehat{P}_{MLE}(10^{-4}) = 0.9995$). In these cases, \textsc{SVDrop}
did not converge to a KKT point. Instead CPAPR oscillated near some other local
minimizer---perhaps a saddle point. \Cref{fig:cpapr-svdrop:fail} demonstrates
that the failure of \textsc{SVDrop} in these two instances was due to setting
the gap ratio tolerance $\gamma$ too large (red $\times$ marker). In both plots, the gap
ratio was always less than the choice of the tolerance in our experiments
($\gamma = 10^6$). When set appropriately, e.g., $\gamma = 8\times 10^3$ (blue
circle marker), Restarted CPAPR with \textsc{SVDrop} converged. Note that the traces
were identical until the gap tolerance was exceeded in the $\gamma
= 8\times 10^3$ case. This triggered a restart, so that the iteration
histories diverged. Observe that the trial in the bottom plot restarted
twice when $\gamma = 8 \times 10^3$. \Cref{alg:cpapr_svdrop} will
always restart until convergence to a KKT point that is not
rank-deficient.

\begin{table}[!ht]
	\centering
	\caption{Convergence results of Restarted CPAPR with \textsc{SVDrop}: the
		number of trials that converged to the MLE, converged to some other KKT
		point, or did not converge to a KKT point. The initial guesses were the
		set of starts from previous experiments where CPAPR without restarting
		(i.e., $\tau=0$) did not converge to the MLE (number of starts
		$N=\text{4,051}$). \emph{Converged} means the solution satisfied the KKT-based
		CPAPR convergence criterion to tolerance $10^{-15}$.}
	\label{tab:cpapr-svdrop:results}
	\resizebox{\textwidth}{!}{
		\begin{tabular}{llrrrrrrrrrrr}
			\toprule
			          &                 & \multicolumn{11}{c}{\textsc{SVDrop}
			inner iterations $\tau$}
			\\ \cmidrule{3-13} Converged & Minimizer       & 0
			& 1    & 2             & 3    & 4 & 5               & 6
			& 7    & 8             & 9 & 10
			\\
			\midrule
			Yes       & MLE             & 0
			& 4024 & \textbf{4049} & 4035 & 4028 & 4029 & 3906 & 3970 & 3983 &
			3990 & 3998 \\
			Yes       & Other KKT point & 3905
			& 0    & 0             & 0    & 0    & 0    & 102  & 43   & 31   &
			24   & 20   \\
			No        & -               & 146
			& 27   & \textbf{2}    & 16   & 23   & 22   & 43   & 38   & 37   &
			37   & 33   \\
			\bottomrule
		\end{tabular}%
	}
\end{table}

\subsubsection{Computational cost}
The formulae to compute computational costs in FLOPS are provided
in~\Cref{sec:costs}. In the best case ($\tau = 2$), the cost of Restarted CPAPR
with \textsc{SVDrop} was 7.04$\times$ higher than CPAPR without restarting
($\tau = 0$).
Although the cost of converging to
the MLE using Restarted CPAPR with \textsc{SVDrop} was more expensive than CPAPR
without restarting, it may be possible to converge to the MLE with higher
probability---$\widehat{P}_{MLE}(10^{-4}) = 0.9995$ versus
$\widehat{P}_{MLE}(10^{-4}) = 0.963$---with the extra work.

\subsubsection{Caveats}
Although our results were promising, we caution that our approach may not always
improve on CPAPR without restarting. There appears to be a ``Goldilocks'' range
for the gap ratio: too large $\gamma$ may decrease the probability that
rank-deficient solutions are identified and too small $\gamma$ may trigger
unwanted restarts.

\paragraph{Choosing $\gamma$ too large}
Starting from the random sample of 10,000 points where CPAPR without restarting
converged to the MLE in previous experiments, \textsc{SVDrop} occasionally
performed worse. When the number of \textsc{SVDrop} inner iterations was
$\tau=1$, 34 trials did not converge to a KKT point. When the number of
\textsc{SVDrop} inner iterations was $\tau=8$, one trial converged to a KKT
point that was not the MLE. In that trial, a restart was triggered when the gap
ratio was greater than $10^6$. However, two gap ratios, which were large
($>2\times 10^4$) but below the threshold, were missed due to the tolerance
having been set too large.

\begin{figure}
	\centering
	\includegraphics[width=\textwidth]{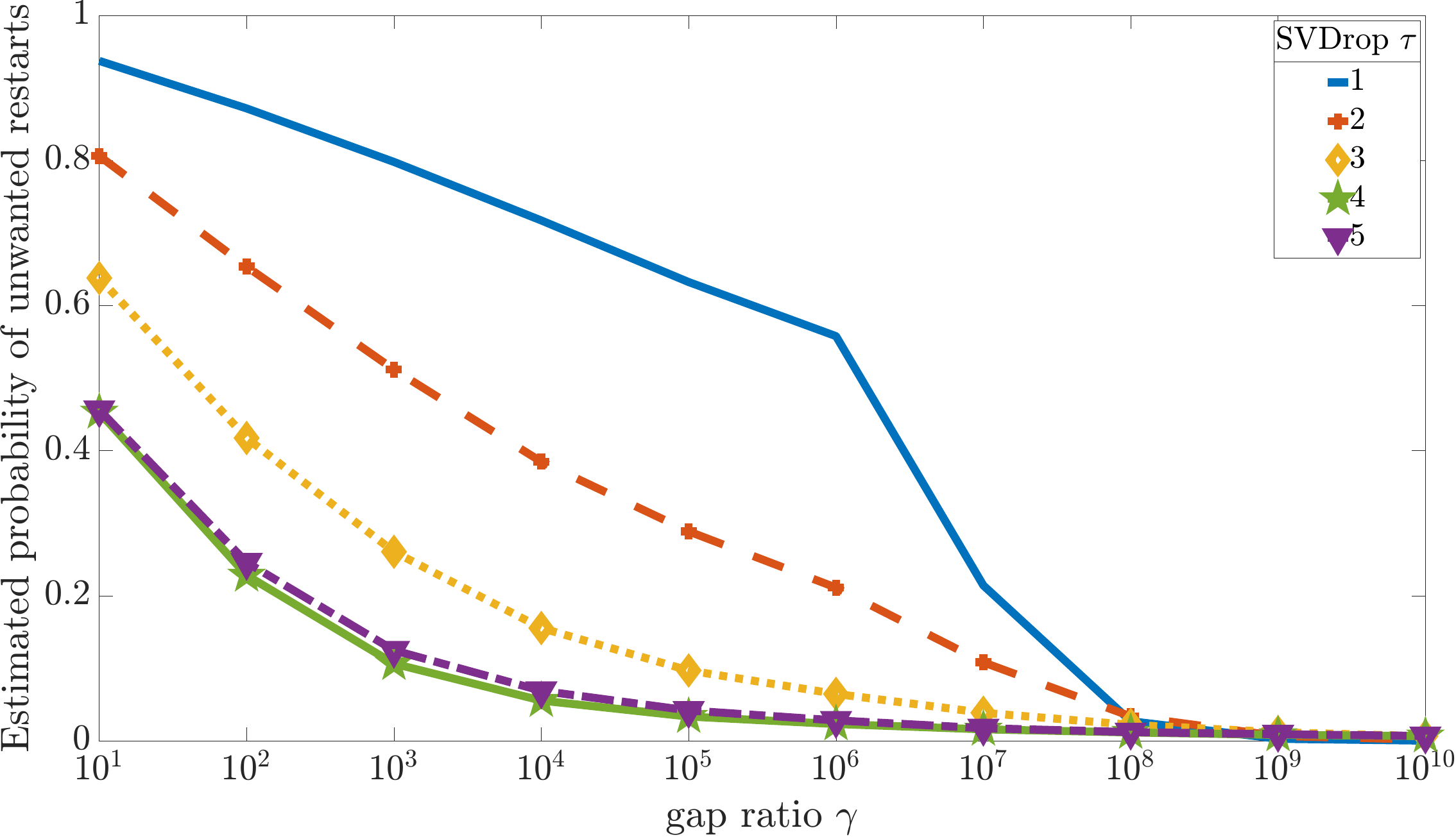}
	\caption{For trials that were known to converge to the MLE without using
	\textsc{SVDrop}, this figure reflects the estimated probabilities that
	\textsc{SVDrop} would trigger an unwanted restart for a range of gap ratios
	$\gamma$.}
	\label{fig:cpapr-svdrop:unwanted-restarts}
\end{figure}

\paragraph{Choosing $\gamma$ too small}
When starting from points known to converge to the MLE without restarting (i.e.,
$\tau = 0$), it is possible that \textsc{SVDrop} could misclassify a solution as
rank-deficient and restart from a new initial guess, which might needlessly
increase the work expended. We consider this unwanted behavior since minimal
computational cost, in addition to accuracy, is a desired characteristic of our
algorithm. \Cref{fig:cpapr-svdrop:unwanted-restarts} shows estimated
probabilities that \textsc{SVDrop} would trigger an unwanted restart for a range
of gap ratios $\gamma$. The implication is that choosing $\gamma$ to be too
small may hurt performance.

\subsection{Discussion}
\label{sec:cpapr-svdrop:discussion}
The results for \textsc{SVDrop} presented here are limited to a small exemplar.
We discuss below several questions that should be addressed before
conclusions about the general efficacy of \textsc{SVDrop} can be considered.

\subsubsection{Connections between $\tau$, $\gamma$, and $l_{max}$} It is possible
that the number of \textsc{SVDrop} steps should be bounded by half the number of
inner iterations, $\tau \leq \lfloor l_{max} / 2 \rfloor$. The reason being that
there is a correspondence between $\tau$, $l_{max}$, and the number of times the
gap ratio is computed per set of inner iterations. To see this, suppose
$l_{max}=10$ and $\tau = 6$. If not converged when $l=6$, then the gap ratio
would be computed only once for that mode; any additional changes to the model
variables beyond that inner iteration would not be captured by \textsc{SVDrop}.
Since rank-deficiency might only be indicated by the singular values of only one
mode, as we observed in \cref{fig:hybrid-gc:experiments:sigmas:traces:mle}, it may
not be apparent that the model had been driven even closer to a rank-deficient
solution while the remaining modes were being optimized. If $\tau=5$, then the
gap ratio would be computed twice: first when $l=5$ and again when $l=10$. In
this case, we would not miss critical changes. Thus it is possible that the gap
ratio should be computed at least twice per set of inner iterations. Two
possible options are to set $\tau = \lfloor l_{max}/2 \rfloor$ or to compute
$\gamma$ when $l = l_{max}$.
We speculate that the choices of $\tau$ and $\gamma$ are highly-coupled when
considering both the probability of convergence to the MLE and the probability
of unwanted restarting. We leave this investigation to future work.

\subsubsection{Swamps}
In future work, we will explore rank-deficiency more systematically for local
minimizers of Poisson CPD problems to better understand the general
applicability of our \textsc{SVDrop} approach. For example,
in~\cref{fig:cpapr-svdrop:inneriters-loop:example}, the CP model appears to
be in a swamp\footnote{The term \emph{swamp} was introduced by Mitchell and
Burdick in~\cite{Mitchell94SlowlyConvergingParafac}. A swamp is \emph{a
phenomenon... in which a [CP] sequence spends a long time in the vicinity of an
inferior resolution before emerging and converging to an acceptable resolution}.
Here we take ``acceptable resolution'' to be the MLE and an ``inferior
resolution'' to be some other local minimizer.} close to some local minimizer
before emerging to converge to the MLE. One direction open to future work is to
compare \textsc{SVDrop} with methods for avoiding 2FD such as those
in~\cite{Giordani16RemediesDegeneracyCandecomp}.

\subsubsection{Efficient computation of gap ratios}
\label{sec:cpapr-svdrop:discussion:svd}
Black-box solvers, e.g., those built on LAPACK
\texttt{xGES*D}\cite{Anderson99LAPACK}, compute \emph{all} of the singular
values of the matrix. \emph{Iterative methods}, such as subspace iteration or
Krylov methods, may be a more practical choice for computing the gap ratio due
to the (relatively) small number of nonzero entries in each unfolding (This is a
consequence of the CPAPR algorithm design that drives elements toward zero.)
Iterative methods are also useful when only a small number of singular values
are needed and are amenable to sparse data structures. While it is less
straightforward to characterize the computational costs and other trade-offs of
iterative methods, it is worth investigating their role in future work.

\section{Conclusions}
\label{sec:conc}

We presented two new methods for Poisson canonical polyadic decomposition:
\textsc{HybridGC} and Restarted CPAPR with \textsc{SVDrop}.

\subsection{\textsc{HybridGC}} Our method can minimize low-rank approximation
error with high accuracy relative to GCP-Adam and CPAPR-MU while reducing
computational costs. Since \textsc{HybridGC} was run	in our experiments with
a far stricter computational budget than GCP-Adam and CPAPR-MU, we argue that
\textsc{HybridGC} can be more computationally efficient. The implication is that
the performance gain allows even more multi-starts, and subsequently, a greater
number of high accuracy approximations. Furthermore, our contribution is a new
method that interpolates two very different algorithms.

\CHANGED{%
One direction for future work is to study the impact on computation
and accuracy of \textsc{HybridGC} when solver tolerances are relaxed.
A second direction is related to the stagnation of \textsc{HybridGC}
in~\cref{fig:hybrid-gc:experiments:nll:data1:traces:loc} near the
local minimizer, which is indicative of
a swamp~\cite{Mitchell94SlowlyConvergingParafac}. We leave it to
future work to compare \textsc{HybridGC} with methods designed to
avoid swamps.
}

\subsection{Restarted CPAPR with \textsc{SVDrop}} Our method identifies
rank-deficient solutions with near-perfect accuracy and has the highest
likelihood of finding the MLE in our experiments. A corollary is that our
algorithm almost always avoids an entire class of minimizers that are different
from the MLE. Our experimental results demonstrate this empirically with
conservative budgets for both restarting and total iteration, alongside other
untuned parameters. Provided more generous allotments and proper tuning, we
expect \textsc{SVDrop} to always identify rank-deficient solutions in the limit
of multi-starts.

\CHANGED{One direction for future work is to extend the CPD ill-conditioning analysis
by Breiding and Vannieuwenhoven~\cite{BrVa18siopt,BrVa18simax,BrVa18amletters}
to the Poisson CPD problem addressed here. Such analysis could prove useful in 
providing a theoretical understanding of the rank-deficient solutions explored empirically in this paper.}

\subsection{Parameter tuning}
Unlike \textsc{SVDrop}, \textsc{HybridGC} lacks a mechanism for changing its
search path. Assuming a pattern behavior exists, it is essential that further
algorithm development uncovers a diagnostic to identify it. Otherwise,
\textsc{HybridGC} will remain reliant on costly ad hoc parameter tuning by the
user. Convergence and the computational cost of Restarted CPAPR with
\textsc{SVDrop} both depend on the complex interplay of search parameters, which
is not well-understood. Although we provided sensible values and rationalized
upper bounds on some parameters, it remains an open question as to how sensitive
\textsc{SVDrop} is to parameter variability. It is essential to better
understand this interplay since \textsc{SVDrop} can be prohibitively expensive
when it does fail. Fortunately, this is rare.

\appendix
\section{\texttt{LowRankSmall} synthetic data tensor and empirical maximum
  likelihood estimator from numerical experiments}
\label{sec:data:low-rank-small}

\renewcommand{\lstlistingname}{Sparse Tensor}
\begin{lstlisting}[
	caption={Synthetic data tensor \texttt{LowRankSmall} used in experiments in~\cref{sec:hybrid-gc:experiments}.},
	label={data:hybrid-gc:low-rank-small-sptensor}
	]
sptensor    %% tensor type
3           % number of dimensions
4 6 8       % sizes of dimensions
17          % number of nonzeros
1 4 1 1     % start of sparse tensor data: <mode-1-index> <mode-2-index> <mode-3-index> <value>
1 4 6 5
1 4 7 9
1 6 6 1
1 6 7 1
2 1 2 6
2 1 4 5
2 1 8 3
2 4 4 1
2 5 2 1
2 6 2 6
2 6 4 8
2 6 8 3
4 1 8 4
4 2 1 1
4 2 8 2
4 5 8 1
\end{lstlisting}

\renewcommand{\lstlistingname}{MLE Kruskal Tensor}
\begin{lstlisting}[
	caption={Kruskal tensor empirical maximum likelihood estimator (MLE) of~\Cref{data:hybrid-gc:low-rank-small-sptensor} computed in experiments in~\cref{sec:hybrid-gc:experiments}. The empirical MLE is the solution with the smallest Poisson loss from among 110,226 random starts.},
	label={data:hybrid-gc:low-rank-small-mle}
	]
ktensor                                                                   %% tensor type
3                                                                         % number of dimensions
4 6 8                                                                     % sizes of dimensions
3                                                                         % number of components
3.2999999999999986e+001 1.7000000000000004e+001 8.0000000000000000e+000   % lambda values
matrix                                                                    %% 1st factor matrix
2                                                                         % number of dimensions
4 3                                                                       % sizes of dimensions
0.0000000000000000e+000 1.0000000000000000e+000 0.0000000000000000e+000   % start of data
1.0000000000000000e+000 0.0000000000000000e+000 7.7243827424339032e-017
0.0000000000000000e+000 0.0000000000000000e+000 0.0000000000000000e+000
9.1691651618349857e-046 0.0000000000000000e+000 1.0000000000000000e+000
matrix                                                                    %% 2nd factor matrix
2                                                                         % number of dimensions
6 3                                                                       % sizes of dimensions
4.2424242424242420e-001 0.0000000000000000e+000 5.0000000000000033e-001   % start of data
9.8952534188649872e-233 0.0000000000000000e+000 3.7499999999999983e-001
0.0000000000000000e+000 0.0000000000000000e+000 0.0000000000000000e+000
3.0303030303030307e-002 8.8235294117647056e-001 0.0000000000000000e+000
3.0303030303030307e-002 0.0000000000000000e+000 1.2499999999999993e-001
5.1515151515151525e-001 1.1764705882352941e-001 8.6248456937576744e-133
matrix                                                                    %% 3rd factor matrix
2                                                                         % number of dimensions
8 3                                                                       % sizes of dimensions
0.0000000000000000e+000 5.8823529411764705e-002 1.2499999999999992e-001   % start of data
3.9393939393939398e-001 0.0000000000000000e+000 3.8294555981328405e-118
0.0000000000000000e+000 0.0000000000000000e+000 0.0000000000000000e+000
4.2424242424242437e-001 0.0000000000000000e+000 1.2073880062159344e-131
0.0000000000000000e+000 0.0000000000000000e+000 0.0000000000000000e+000
0.0000000000000000e+000 3.5294117647058826e-001 0.0000000000000000e+000
0.0000000000000000e+000 5.8823529411764708e-001 0.0000000000000000e+000
1.8181818181818168e-001 0.0000000000000000e+000 8.7500000000000011e-001
\end{lstlisting}

\section{Computational cost derivations}
\label{sec:costs}

\subsection{CPAPR-MU operation count}
CPAPR-MU was the first algorithm developed for Poisson CPD and it is important
to assess its cost in terms of the number of floating point operations (FLOPS)
required. We are only interested in sparse tensor computations, so we do not
consider the costs of operations for dense tensors.

The work in an iteration of CPAPR-MU is dominated by the following operations:
\begin{enumerate}
	\item	Sequence of Khatri-Rao products: $\mat{\Pi}   \gets ( \mat{A}^{(d)}
		      \odot \dots \odot \mat{A}^{(n+1)} \odot \mat{A}^{(n-1)} \odot
		      \dots \odot \mat{A}^{(1)} )^T$
	\item  Implicit MTTKRP: $\mat{\Phi}  \gets \left( \mat{X}_{(k)} \oslash
		      \left(\mat{B}\mat{\Pi}\right)\right)\mat{\Pi}^T$
	\item Multiplicative update: $\mat{B}     \gets \mat{B} \ast \mat{\Phi}$
\end{enumerate}

The first line is a sequence of Khatri-Rao products, where the binary operator
$\odot$ denotes the Khatri-Rao product between matrices. The Khatri-Rao product
is a primary component of the \emph{matricized tensor times Khatri-Rao product}
(MTTKRP), a key computational kernel in many tensor algorithms, not just CPD.
Improving performance of the MTTKRP is a very active research area. For our
purposes, we treat Khatri-Rao product as a black box and do not count its costs.

The second line computes the $\mat{\Phi}$ matrix used in the multiplicative
update. Note that $\mat{B}$ is simply the $n$-th factor matrix $\mat{A}^{(n)}$
scaled by the $\lambda$ weights, i.e., $\mat{B} =
	\mat{A}^{(n)}\diag(\mat{\lambda})$. We write that it is an implicit MTTKRP since
1) and 2) together are mathematically equivalent to MTTKRP. The MTTKRP is
efficient because it can be done without forming dense arrays and the implicit
MTTKRP can be performed even more efficiently due to the special structure of
the matricized tensor in minimizing the Poisson loss for sparse tensors. The
matrix multiplication $\mat{B}\mat{\Pi}$ requires $\BigO(R \prod_{k=1}^{d} I_k)$
arithmetic. The elementwise division $\mat{X}_{(k)} \oslash
	\left(\mat{B}\mat{\Pi}\right)$ depends on the number of nonzeros in $\tns{X}$;
thus it requires $\nnz(\tns{X})$ operations. Lastly, the product $\left(
	\mat{X}_{(k)} \oslash \left(\mat{B}\mat{\Pi}\right)\right)\mat{\Pi}^T$ requires
$\BigO(R \prod_{k=1}^{d} I_k)$ arithmetic.

The elementwise multiplication in the third line, $\mat{B} \ast \mat{\Phi}$,
requires $RI_n$ multiplications in the $n$-th mode. All together, the number of
FLOPS required per inner iteration for the $n$-th mode is
\begin{equation}
	\sim \nnz(\tns{X}) + rI_n + 2R \prod_{k=1}^{d} I_k \text{ FLOPS}.
\end{equation}

Note: we ignored the following computations with negligible costs (with
corresponding line numbers from~\cite[Alg.
	3]{Chi12TensorsSparsityNonnegative}):
\begin{enumerate}
	\item inadmissible zero avoidance (Line 4);
	\item the shift of weights from $\mat{\Lambda}$ to mode-$n$ and vice versa
	      (Lines 5, 15, 16); and,
	\item the check for convergence (Line 9).
\end{enumerate}

\section{Cyclic GCP-CPAPR}
\label{sec:cyclic-gc}

We develop Cyclic GCP-CPAPR (\textsc{CyclicGC}), a generalized form of
\textsc{HybridGC} that \emph{cycles} between a stochastic method to compute a
model approximation and a deterministic method to resolve the model to the best
accuracy possible at scale. In our formulation, \textsc{HybridGC} is
\textsc{CyclicGC} with a single cycle. We define parameterizations and cycle
\emph{strategies}, which prescribe how \textsc{CyclicGC} iterates in each cycle.

Let $L \in \N$ be a number of \emph{cycles}. Define \emph{strategy} to be the
$L$-length array of structures, \texttt{strat}, specifying the following for
each cycle $l \in \{ 1, \ldots, L \}$:
\begin{itemize}
	\item \texttt{S\_opts}: stochastic search parameterization, including solver
	      and search budget, $j$, measured in outer iterations.
	\item \texttt{D\_opts}: deterministic search parameterization, including
	      solver and search budget, $k$, measured in outer iterations.
\end{itemize}
\textsc{CyclicGC} iterates from an initial guess $\tns{M}^{(0)}$ via a two-stage
alternation between stochastic and deterministic search for $L$ cycles to return
a Poisson CP tensor approximation $\tns{\widehat{M}}$ as an estimate to
$\tns{M}^*$. In the first stage of the $l$-th cycle, the stochastic solver
iterates from $\tns{M}^{(l-1)}$ for $j$ outer iterations, parameterized by
\texttt{strat(l).S\_opts} to return an intermediate solution, $\tns{M}^{(l)}$.
In the second stage, the deterministic solver refines $\tns{M}^{(l)}$ for $k$
outer iterations, parameterized by \texttt{strat(l).D\_opts}, to return the
$l$-th iterate, $\tns{M}^{(l)}$, overwriting the output from the previous stage.

\begin{algorithm}
	\caption{Cyclic GCP-CPAPR}
	\label{alg:cyclic-gc}
	\begin{algorithmic}[1]
		\Function{CyclicGC}{tensor $\tns{X}$, rank $r$, initial guess
			$\tns{M}^{(0)}$, number of cycles $L$, $L$-array of structures
			\texttt{strat} defining $L$ strategies.} \For{$l=1,\ldots,L$} \State
			$\tns{M}^{(l)} \gets \textsc{GCP}(\tns{X}, \,r, \,\tns{M}^{(l-1)}\,
			\texttt{strat(l).S\_opts})$ \State $\tns{M}^{(l)} \gets
			\textsc{CPAPR}(\tns{X}, \,r, \,\tns{M}^{(l)}, \,
			\texttt{strat(l).D\_opts})$ \EndFor \State \Return model tensor
			$\tns{\widehat{M}} = \tns{M}^{(L)}$ as estimate to $\tns{M}^*$
			\EndFunction
	\end{algorithmic}
\end{algorithm}

\section{Supplemental numerical results}
\label{sec:hybrid-gc:experiments:supp}

This section presents additional numerical results that are supplementary to
those in \cref{sec:hybrid-gc:experiments}.

\begin{figure}[!ht]
    \centering
    \begin{subfigure}[t]{0.32\textwidth}
        \includegraphics[width=\textwidth]{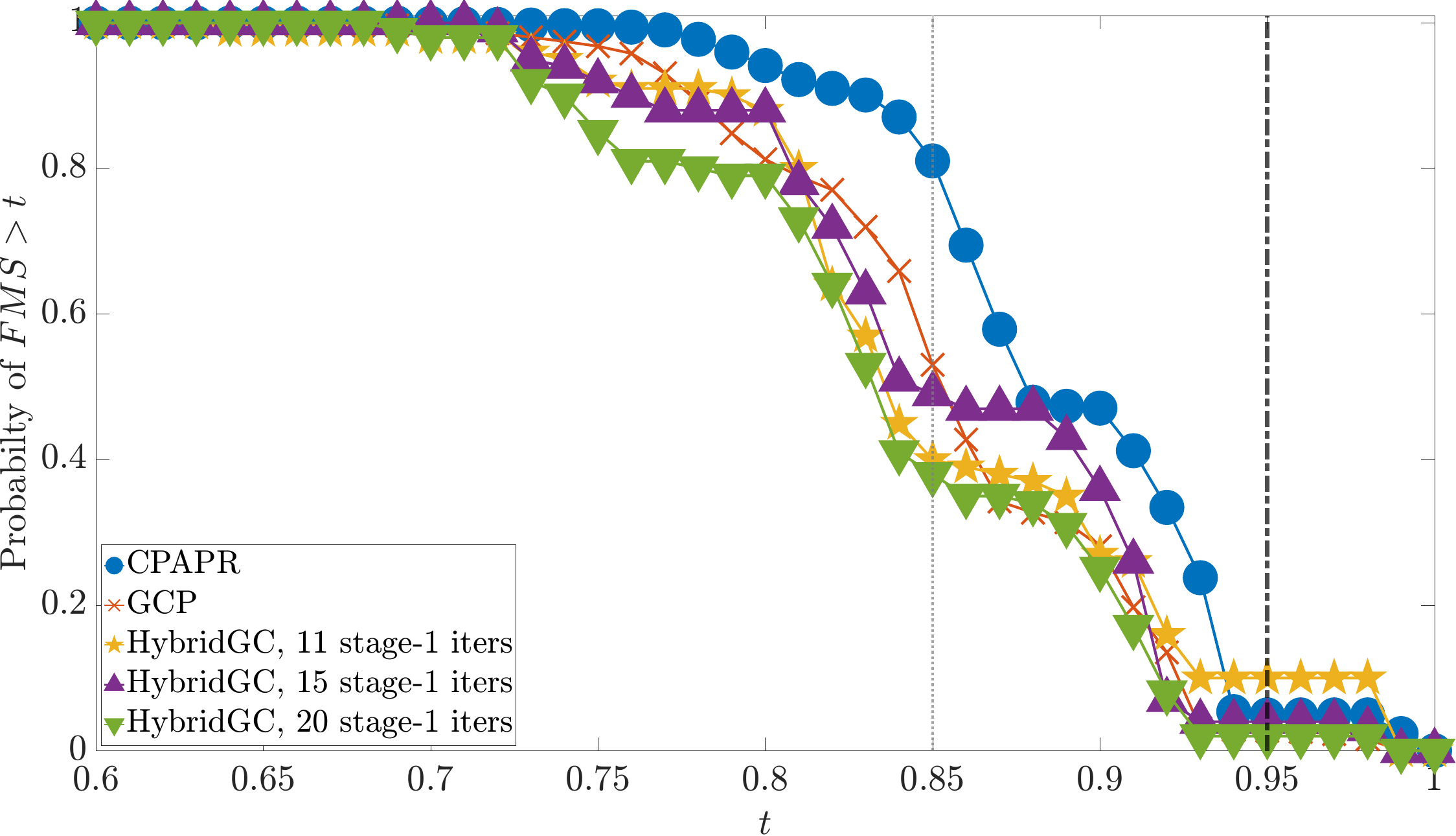}
        \label{fig:hybrid-gc:experiments:fms:prob_mle_est:2}  \caption{11 to
            20.} \end{subfigure}\hfill
    \begin{subfigure}[t]{0.32\textwidth}
        \includegraphics[width=\textwidth]{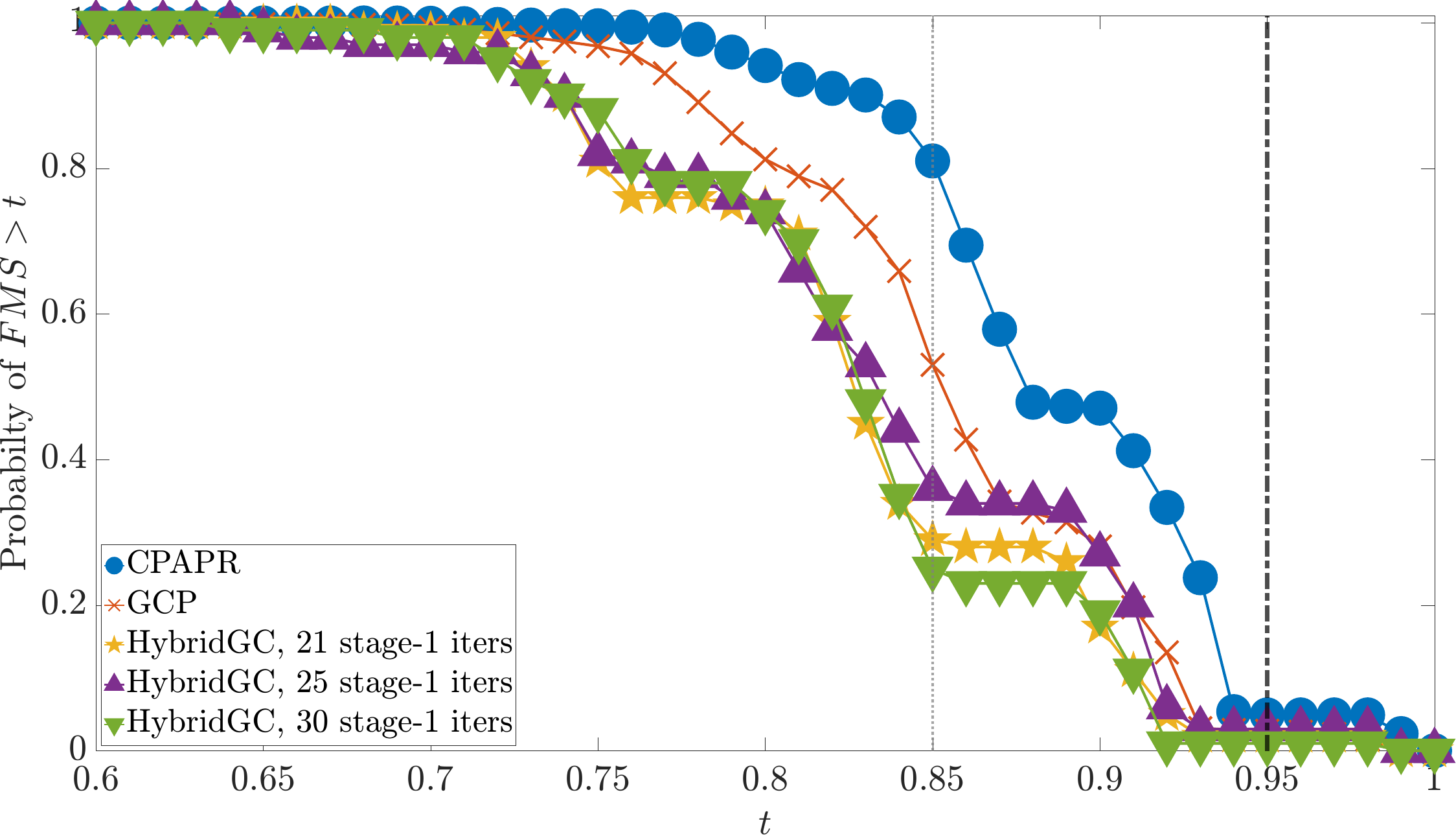}
        \label{fig:hybrid-gc:experiments:fms:prob_mle_est:3}  \caption{21 to
            30.} \end{subfigure}\hfill
    \begin{subfigure}[t]{0.32\textwidth}
        \includegraphics[width=\textwidth]{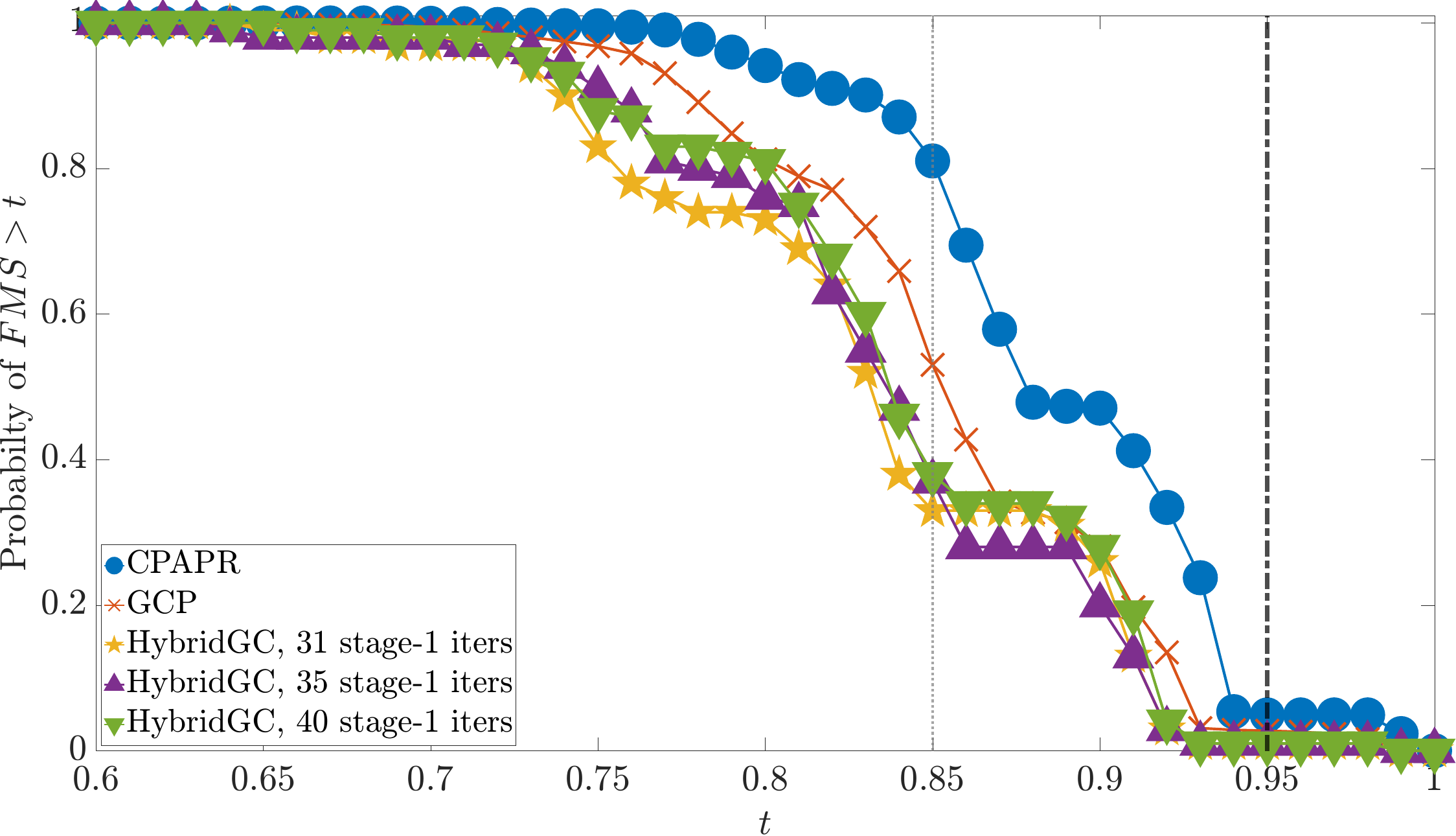}
        \label{fig:hybrid-gc:experiments:fms:prob_mle_est:4}  \caption{31 to
            40.} \end{subfigure} \\
    \begin{subfigure}[t]{0.32\textwidth}
        \includegraphics[width=\textwidth]{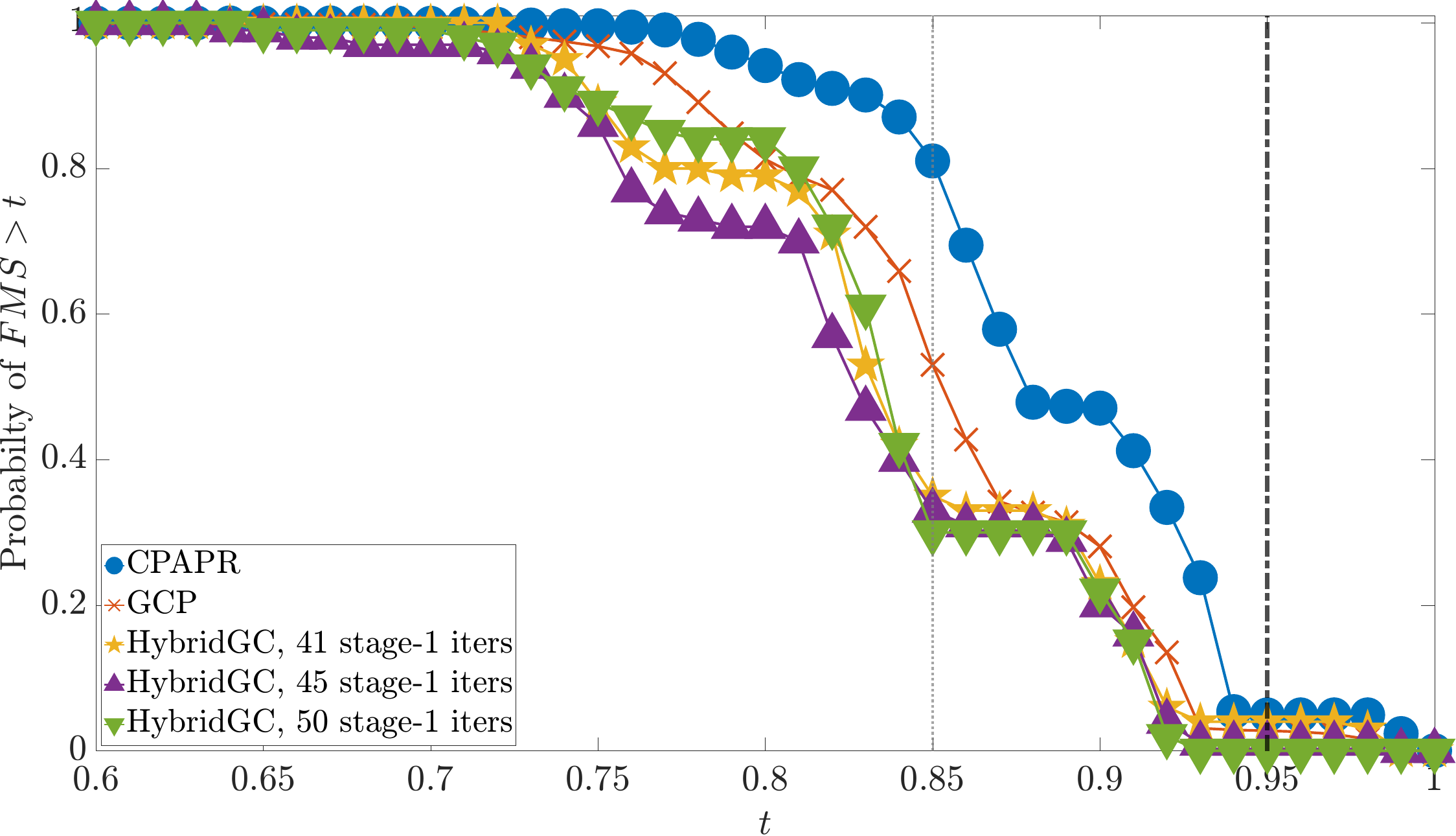}
        \label{fig:hybrid-gc:experiments:fms:prob_mle_est:5}  \caption{41 to
            50.} \end{subfigure}\hfill
    \begin{subfigure}[t]{0.32\textwidth}
        \includegraphics[width=\textwidth]{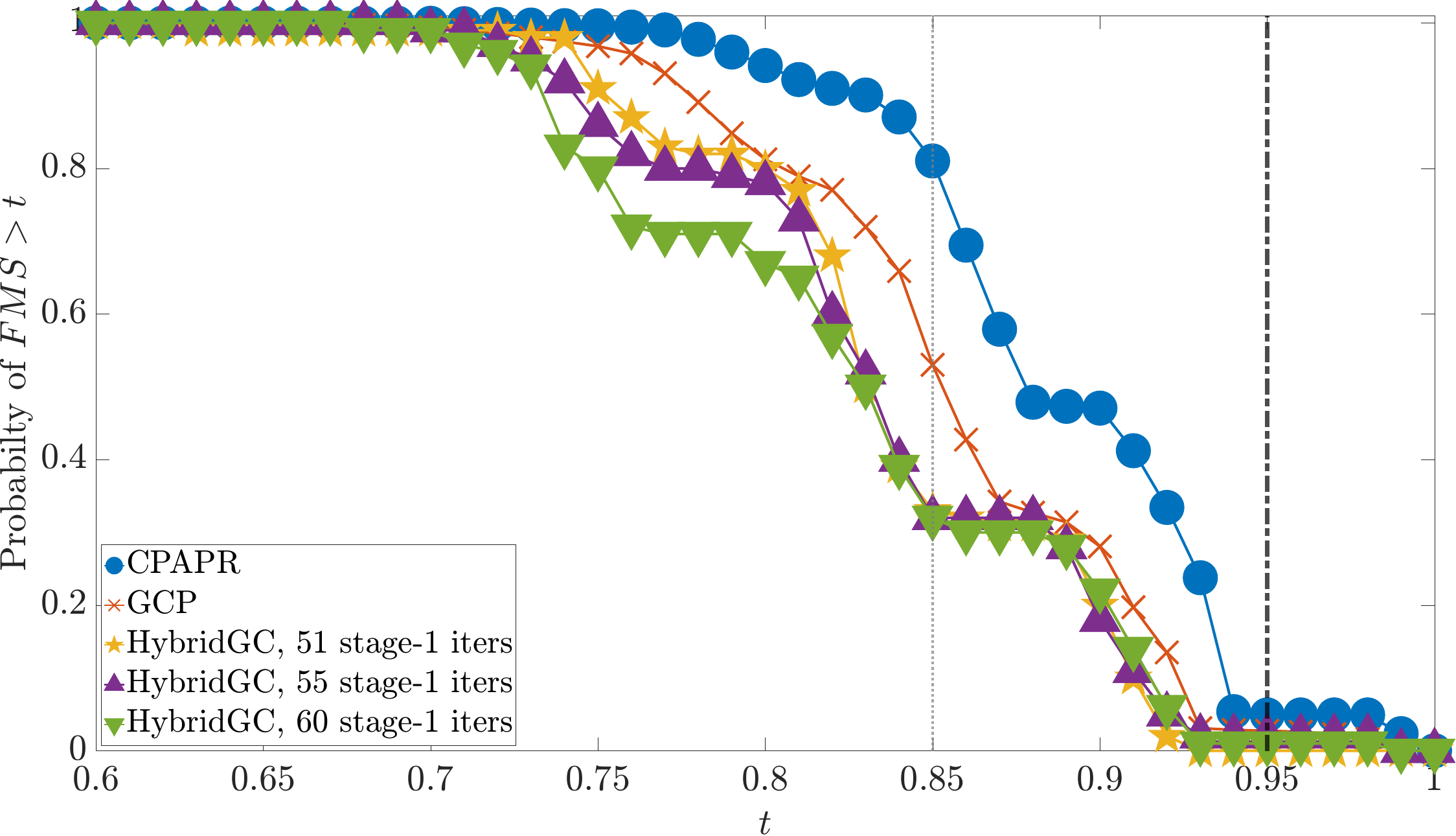}
        \label{fig:hybrid-gc:experiments:fms:prob_mle_est:6}  \caption{51 to
            60.} \end{subfigure}\hfill
    \begin{subfigure}[t]{0.32\textwidth}
        \includegraphics[width=\textwidth]{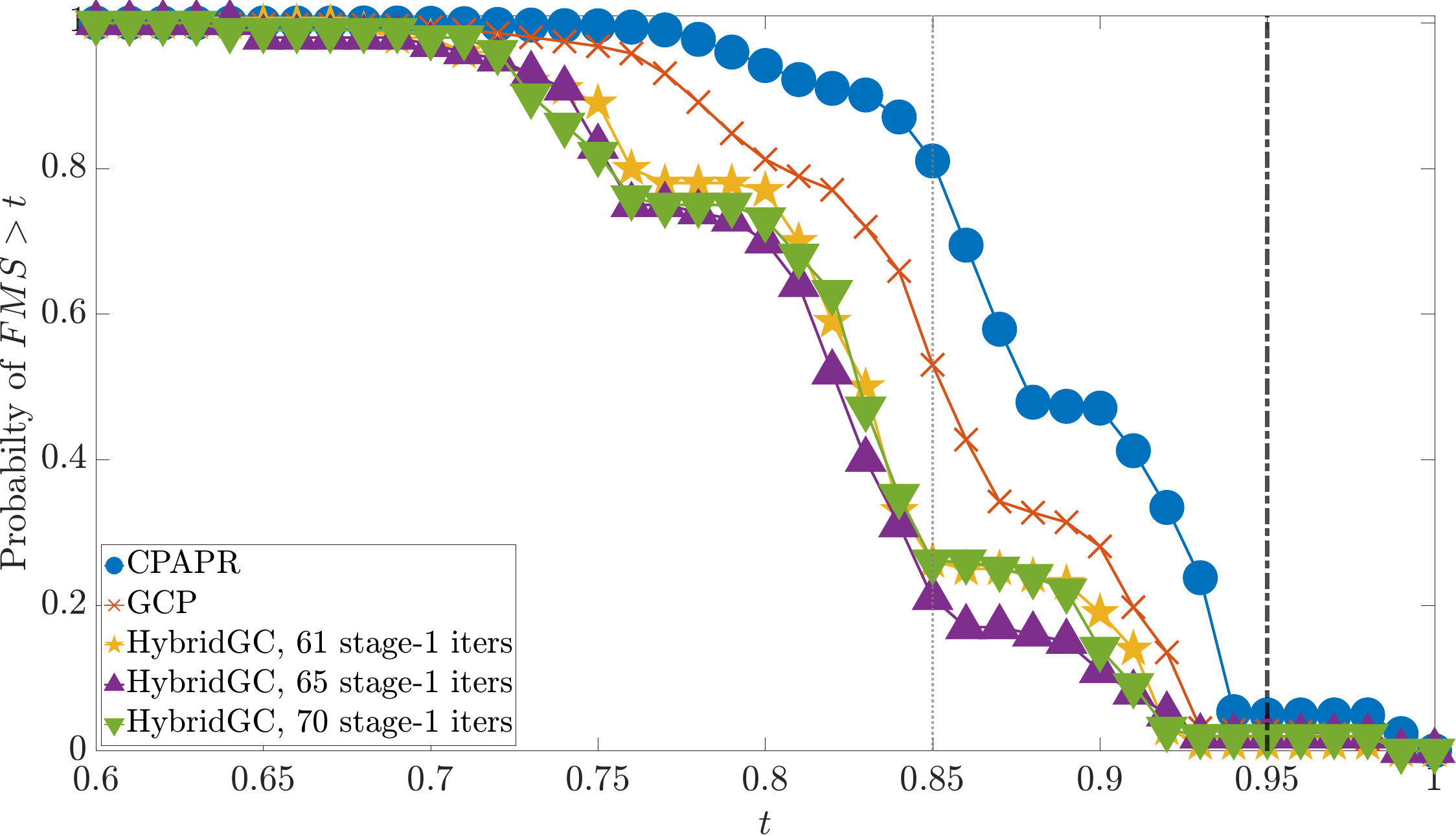}
        \label{fig:hybrid-gc:experiments:fms:prob_mle_est:7}  \caption{61 to
            70.} \end{subfigure} \\
    \begin{subfigure}[t]{0.32\textwidth}
        \includegraphics[width=\textwidth]{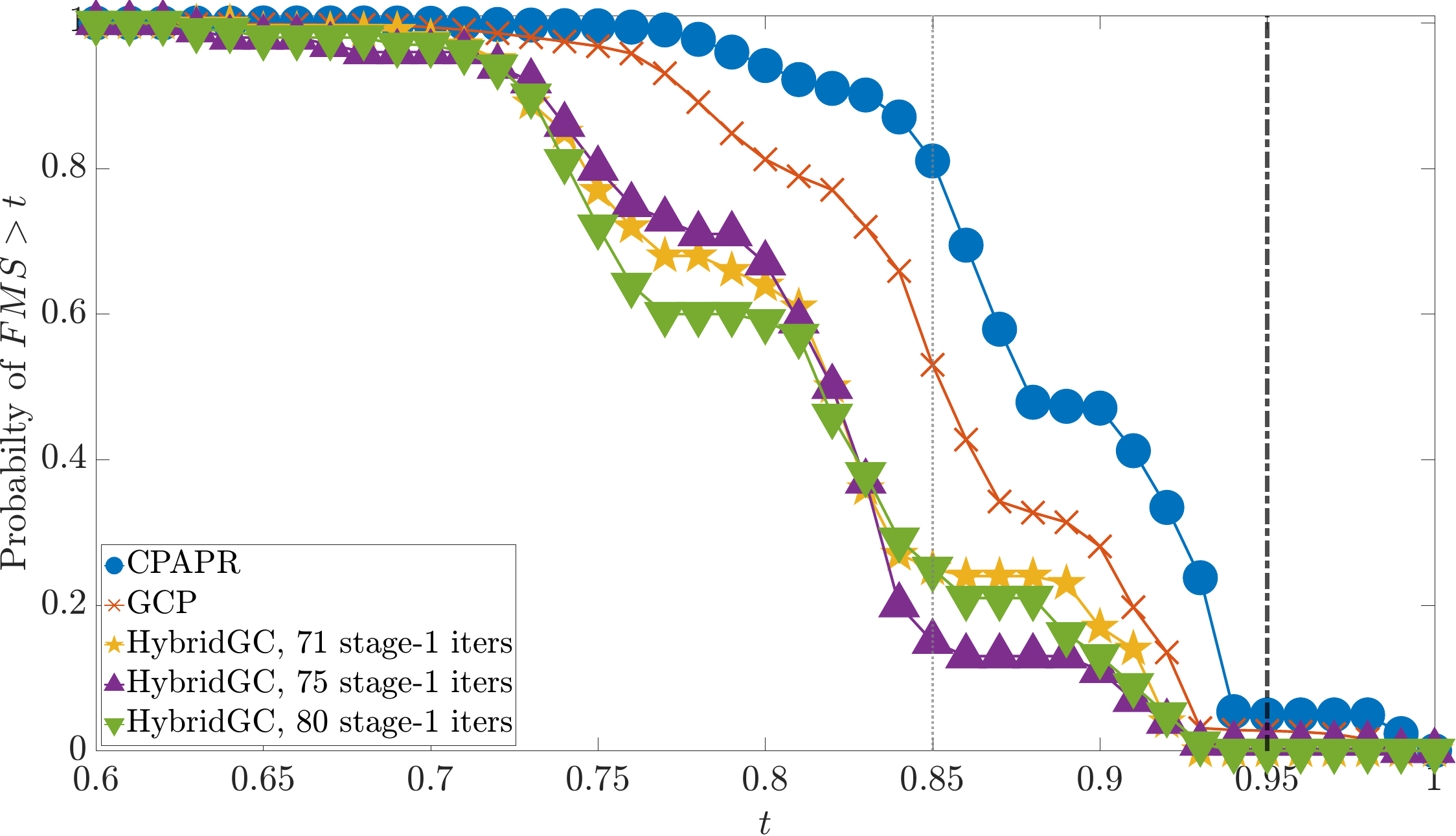}
        \label{fig:hybrid-gc:experiments:fms:prob_mle_est:8}  \caption{71 to
            80.} \end{subfigure}\hfill
    \begin{subfigure}[t]{0.32\textwidth}
        \includegraphics[width=\textwidth]{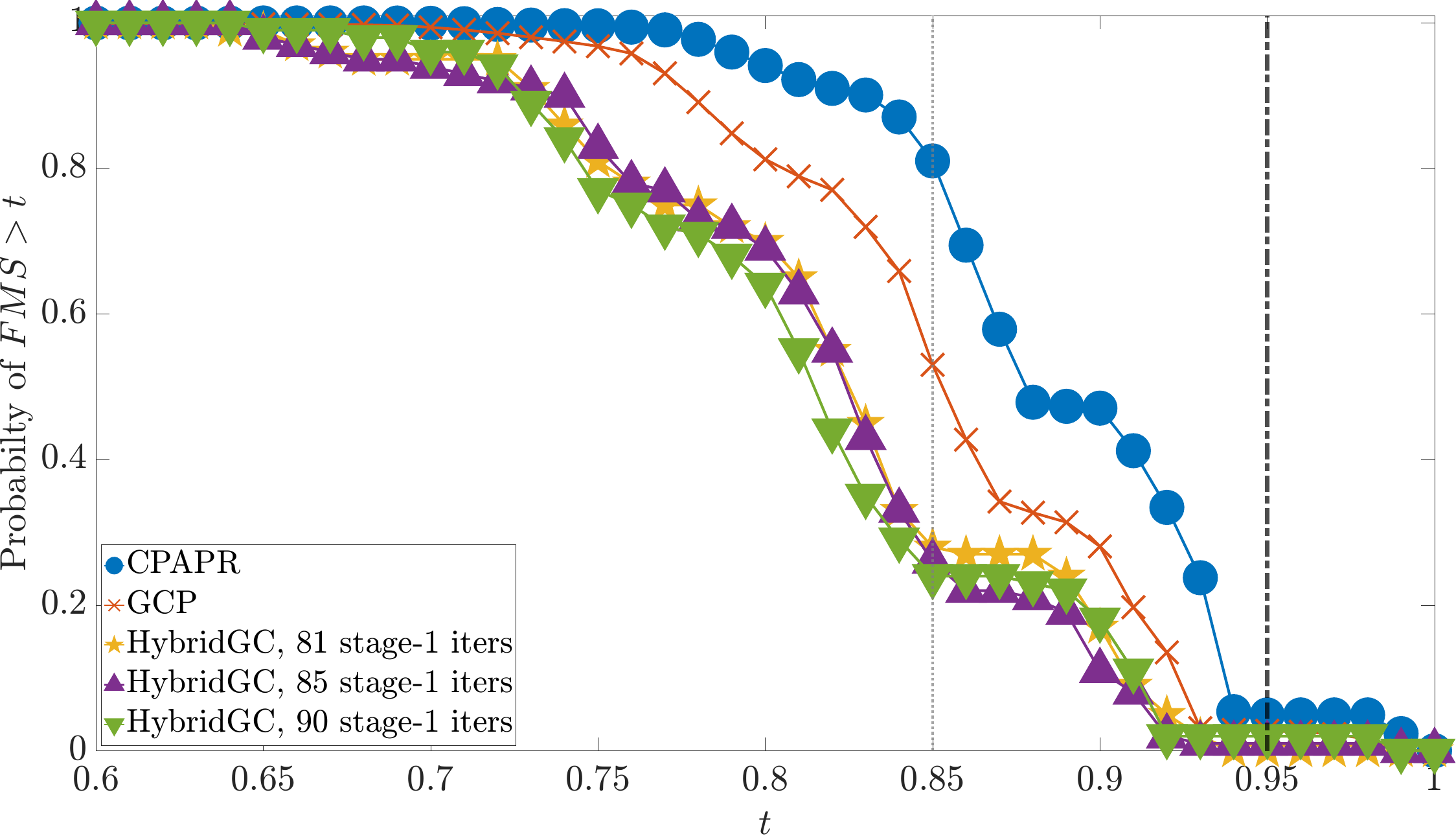}
        \label{fig:hybrid-gc:experiments:fms:prob_mle_est:9}  \caption{81 to
            90.} \end{subfigure}\hfill
    \begin{subfigure}[t]{0.32\textwidth}
        \includegraphics[width=\textwidth]{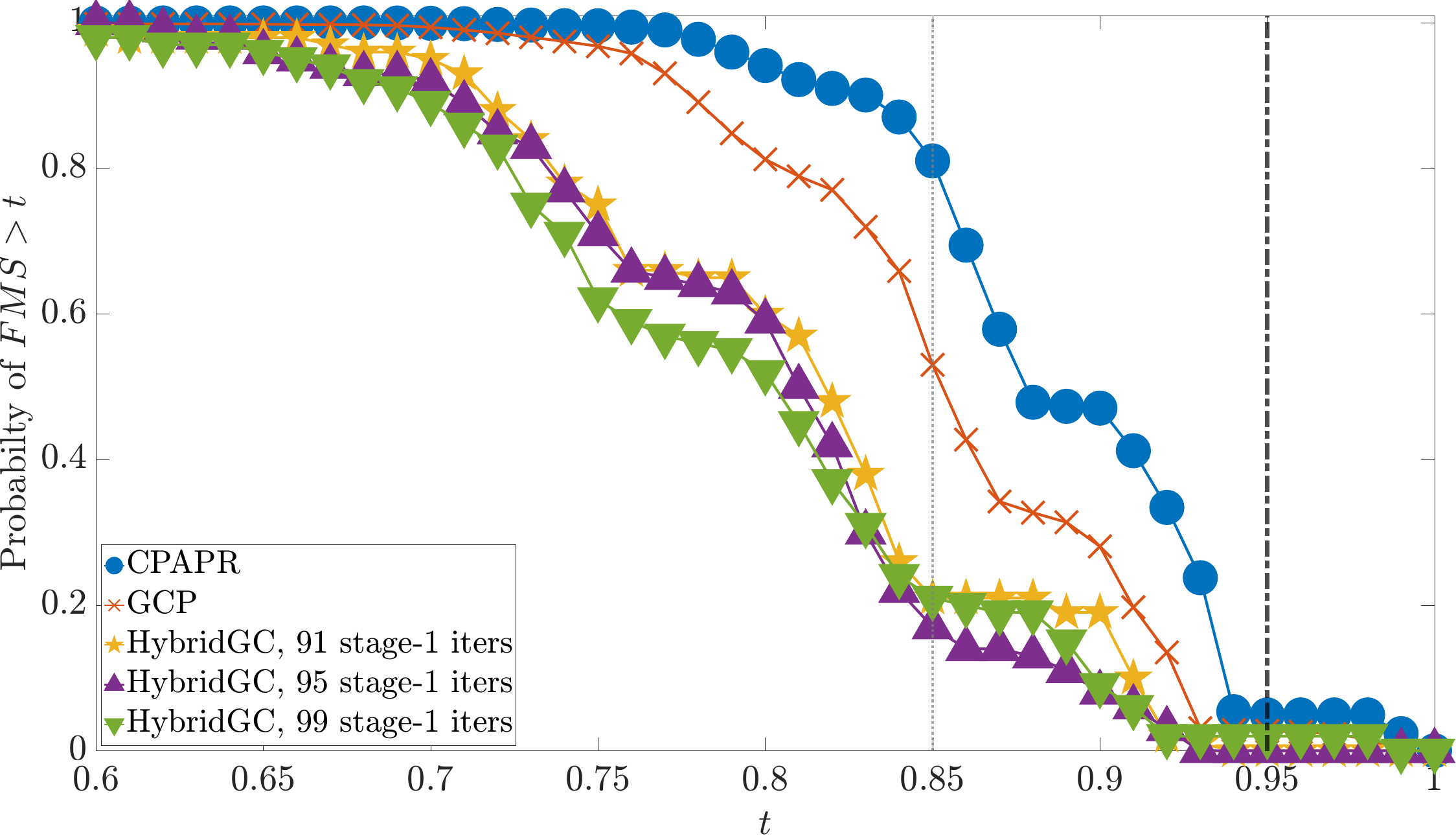}
        \label{fig:hybrid-gc:experiments:fms:prob_mle_est:10} \caption{91 to
            99.}
    \end{subfigure}
    \caption{Factor match scores between CP models computed with HybridGC,
        CPAPR-MU, and GCP-Adam and the approximate global optimizer,
        $\tns{\widehat{M}}_{\mathcal{S}}^*$. The dash-dot gray vertical lines
        and dotted black vertical lines denote the levels of ``similar'' and
        ``equal'' described in~\cite{Lorenzo-Seva06TuckerCongruenceCoefficient}.
        Colormaps scaled for clarity.}
    \label{fig:hybrid-gc:experiments:fms:prob_mle_est:2-10}
\end{figure}

\clearpage
\section*{Acknowledgments}
\CHANGED{%
We thank Eric Phipps of Sandia National Laboratories for assistance
with Genten, a high-performance GCP solver. 
Sandia National Laboratories is a multi-mission laboratory managed and operated
by National Technology \& Engineering Solutions of Sandia, LLC (NTESS), a wholly
owned subsidiary of Honeywell International Inc., for the U.S. Department of
Energy’s National Nuclear Security Administration (DOE/NNSA) under contract
DE-NA0003525. This written work is authored by an employee of NTESS. The
employee, not NTESS, owns the right, title and interest in and to the written
work and is responsible for its contents. Any subjective views or opinions that
might be expressed in the written work do not necessarily represent the views of
the U.S. Government. The publisher acknowledges that the U.S. Government retains
a non-exclusive, paid-up, irrevocable, world-wide license to publish or
reproduce the published form of this written work or allow others to do so, for
U.S. Government purposes. The DOE will provide public access to results of
federally sponsored research in accordance with the DOE Public Access Plan.
}

\clearpage
\bibliographystyle{siamplain}
\bibliography{refs}

\end{document}